\documentclass[times]{nmeauthmodi}

\usepackage{algorithm}
\usepackage{algorithmic}
\usepackage{amsmath,amssymb,amsfonts,amsthm}
\usepackage{mathrsfs}
\usepackage{bm}
\usepackage[overload]{empheq}
\usepackage{epstopdf}
\usepackage{moreverb,bm,psfrag,color} 
\usepackage{subfigure}
\usepackage{color}
\usepackage{pstcol} 

\graphicspath {{./}{./Figures/}}

\newtheorem{theorem}{Theorem}[section]
\def\volumeyear{}

\begin{document}

\runningheads{L.~Fedeli, M.~Ortiz and A. Pandolfi}{Optimal transportation theory and advection-diffusion problems}

\title{Geometrically-exact time-integration mesh-free schemes for advection-diffusion problems derived from optimal transportation theory and their connection with particle methods}

\author{L.~Fedeli\affil{1}, A.~Pandolfi\affil{2}, M.~Ortiz\affil{3}}

\address{
\affilnum{1}{University of Pennsylvania,
Mechanical Engineering and Applied Sciences Department, Philadelphia PA, 19104, USA}
\\
\affilnum{2}{Politecnico di Milano,
Civil and Environmental Engineering Department, 20133 Milano, Italy}
\\
\affilnum{3}
{California Institute of Technology,
Engineering and Applied Science Division, Pasadena CA, 91125, USA}
}

\corraddr{anna.pandolfi@polimi.it}

\begin{abstract}
We develop an Optimal Transportation Meshfree (OTM) particle method for advection-diffusion in which the concentration or density of the diffusive species is approximated by Dirac measures. We resort to an incremental variational principle for purposes of time discretization of the diffusive step. This principle characterizes the evolution of the density as a competition between the Wasserstein distance between two consecutive densities and entropy. Exploiting the structure of the Euler-Lagrange equations, we approximate the density as a collection of Diracs. The interpolation of the incremental transport map is effected through mesh-free max-ent interpolation. Remarkably, the resulting update is geometrically exact with respect to advection and volume. We present three-dimensional examples of application that illustrate the scope and robustness of the method.
\end{abstract}


\maketitle

\section{Introduction}

Advection-diffusion problems have been traditionally regarded as defining time-evolutions in linear spaces. For instance, the classical analysis of the heat equation regards its solutions as elements of the Sobolev space $L^2(0,T;H^1_0(\Omega))$ with time derivatives in $L^2(0,T;H^{-1}(\Omega))$, \cite{Evans1998}. The traditional formulation of numerical approximation schemes for advection-diffusion problems correspondingly works fully within this linear space structure and aims to approximate the governing PDEs by recourse to, e.~g., finite-difference or finite-element discretization schemes. Such schemes are inevitably beset by a toxic combination of spurious oscillations, instabilities, poor accuracy, lack of monotonicity, numerical diffusion, leaking of mass and others (e.~g., \cite{Benson1992}). Despite massive and sustained efforts to overcome them, such deficiencies have shown remarkably stubborn staying power.

More recently, it has been recognized, starting with the seminal work of Jordan, Kinderlehrer and Otto \cite{JordanKinderlehrerOtto1997, JordanKinderlehrerOtto1998, JordanKinderlehrerOtto1999}, that the linear functional-analytical framework is not-at-all natural or convenient and that a more natural and productive framework is to regard advection-diffusion problems as transport problems defining time-evolutions in spaces of measures. From this new perspective, the essential difficulties experienced by linear-space-based approximation schemes for advection-diffusion problems maybe be regarded as the 'punishment' for the 'crime' of formulating the problems in an unnatural linear-space setting. Indeed, as we shall see, the overwhelming advantage of transport methods is that they are {\sl geometrically exact}, i.~e., they are exact with respect to advection and volume. This property effectively eliminates in one fell swoop the difficulties experienced by traditional linear-space methods in dealing with advection and volume, albeit at the expense of working in a somewhat more challenging non-linear space, or manifold, framework.

In a parallel development, particle methods have attained considerable acceptance in solid and fluid mechanics, e.~g, the Smoothed Particle Hydrodynamics method (SPH) \cite{lucy:1977, monaghan:1992}, the Material Point Method \cite{Sulsky94}, the Reproducing Kernel Particle method \cite{liu:1995}, the Corrective Smoothed Particle Method \cite{chen:1999}, the Modified Smoothed Particle Hydrodynamics method \cite{zhang:2004}, the Optimal-Transportation Meshfree method (OTM) \cite{LiHabbalOrtiz2010}, and others. Particle methods have the appeal of reducing problems to the evolution of interacting discrete particles, in seeming analogy with molecular dynamics. There is also a clear connection between particle methods and transport of measures. Indeed, the particles maybe regarded as Dirac-delta approximations of otherwise continuous measures, or densities. However, continuum transport problems differ from true particle-dynamics problems in one important respect: in continuum transport problems, the velocity field that governs the instantaneous motion of the particles depends on local {\sl density gradients}, not just particle positions. This differential structure introduces additional regularity requirements of the density interpolation, such as conformity, and is at the core of the difficulties that plague particle methods such as SPH. The differential structure of the mobility law of continuum transport problems in fact necessitates two types of representations: one for the measure itself, e.~g., represented as a collection of Diracs; and a more regular representation for the transport velocity field, e.~g., based on meshfree conforming interpolation. This double representational requirement was recognized in connection with the OTM method \cite{LiHabbalOrtiz2010}, but has remained obscure in much of the literature on particle methods.

The objective of the present work is to develop an OTM particle method for advection-diffusion in which the concentration or density of the diffusive species is approximated by Dirac measures. The new OTM formulation hybridizes elements of a Galerkin approximation with those of an updated Lagrangian approach in the context of optimal transport theory, to provide an alternative to traditional schemes formulated in linear spaces. We resort to the JKO variational principle \cite{JordanKinderlehrerOtto1997, JordanKinderlehrerOtto1998, JordanKinderlehrerOtto1999} for purposes of time discretization of the diffusive step. This principle characterizes the evolution of the density as a competition between the Wasserstein distance, which penalizes departures from the initial conditions, and entropy, which tends to spread the density and make it uniform over the domain. Remarkably, the resulting update is geometrically exact with respect to advection and volume. We specifically regard the JKO incremental functional as a functional of an incremental transport map which rearranges the density over the time step. Finally, we proceed to discretize the JKO functional in space. This step requires two types of discretization, one for the density and another for the incremental transport map. Exploiting the structure of the Euler-Lagrange equations, which are linear in the density, we may treat said density as a measure and approximate it as a collection of Diracs. The interpolation of the transport map requires more regularity as, in particular, the transport map carries Jacobian information. We specifically resort to the max-ent interpolation scheme proposed by \cite{ArroyoOrtiz2006}. This interpolation scheme is meshfree and conforming and supplies converging approximations in general $W^{1,p}$ spaces \cite{Bompadre2012a}.

The paper is organized as follows. In Sec.~\ref{sec:TimeDiscretization}, we describe how to reformulate the classical diffusion equations as an optimal transportation problem, specifically how a variational structure naturally emerges from time discretization. In Sec.~\ref{sec:SpaceDiscretization}, we present the approximation used to discretize the density and the incremental transport map. In Sec.~\ref{sec:3DExamples} we present numerical schemes that illustrate the scope and properties of the resulting OTM approach.  Finally, in the Appendix we collect the main mathematical elements that {\black underpin the} present approach.

\section{Formulation of the method of approximation}
\label{sec:TimeDiscretization}

We consider the advection-diffusion initial-boundary-value problem
\begin{subequations}\label{system_complete}
\begin{align}[left = \empheqlbrace\,]
    &
    \partial_t \rho + \nabla\cdot(\rho u) = \kappa \Delta \rho,
    & \text{in } \Omega \times [0,T] ,
    \label{eq:TD:Diff1}
    \\ &
    \kappa \nabla\rho \cdot n = 0,
    & \text{in } \partial\Omega \times [0,T] ,
    \label{eq:TD:Diff2}
    \\ &
    \rho(x,0) = \rho_0(x),
    & \text{in } \Omega ,
    \label{eq:TD:Diff3}
\end{align}
\end{subequations}
where $\rho$ is the unknown density or concentration, $\rho_0$ is its initial value, $u$ is a given advection velocity field, $\kappa$ is a diffusion coefficient, $\Omega$ is a bounded domain in $\mathbb{R}^d$ and $n$ is the outward unit normal at the boundary.

As stated, problem (\ref{system_complete}) is defined for densities of sufficient regularity in an appropriate linear space (cf., e.~g., \cite{Evans1998}). In particular, {\black equation \eqref{system_complete}} does not make sense for general measures such as Dirac masses and, therefore, it is not in a form suitable for the formulation of particle methods.

In order to eliminate this obstacle, we proceed to reformulate problem (\ref{system_complete}) as a problem of transport of measures.
We begin by noting that (\ref{eq:TD:Diff1}) can be equivalently reformulated as the pair of first-order partial-differential equations
\begin{subequations}\label{eq:TD:RV}
\begin{align}
    & \label{eq:TD:RV1}
    \partial_t \rho + \nabla\cdot(\rho v) = 0 ,
    \\ & \label{rVre5S}
    \rho v = \rho u - \kappa \nabla\rho ,
\end{align}
\end{subequations}
where $v$ is a velocity field that results from the combined effect of advection and diffusion. We identify (\ref{eq:TD:RV1}) as the continuity equation of Eulerian continuum mechanics and (\ref{rVre5S}) as a mobility law combining the effects of advection and diffusion. The corresponding Lagrangian formulation of the transport problem is
\begin{subequations}
\begin{align}
    & \label{2uN5WW}
    \rho(\varphi(x,t),t)
    =
    \frac{\rho_0(x)}{\det\nabla\varphi(x,t)}
    \\ & \label{yC2krR}
    \partial_t \varphi(x,t)
    =
    v(\varphi(x,t),t) ,
\end{align}
\end{subequations}
where $\varphi : \Omega \times [0,T] \to \Omega$ is the transport map. Eq.~(\ref{2uN5WW}) simply states that $\rho(x,t)$ is the pushforward of $\rho_0(x)$ by the transport map $\varphi(\cdot,t)$, whereas eq.~(\ref{yC2krR}) relates the transport map to the particle velocity.

As stated, the transport problem (\ref{eq:TD:RV}) still requires regularity of the density, but a further weakening of the equations permits extending them to general Radon measures. Thus, (\ref{2uN5WW}) can be reformulated more generally in integral form as
\begin{equation}\label{dZ9cNP}
    \int_\Omega \eta(y) \rho(y,t) \, dy
    =
    \int_\Omega \eta(\varphi(x,t)) \rho_0(x) \, dx ,
\end{equation}
for all continuous test functions $\eta$ with compact support in $\Omega$.
We note that the densities $\rho_0$ and $\rho$ now enter linearly into the integrals. Provided that $\varphi$ has sufficient regularity, (\ref{dZ9cNP}) admits the extension to measures
\begin{equation}\label{UjA9vT}
    \int_\Omega \eta(y) d\mu(y,t)
    =
    \int_\Omega \eta(\varphi(x,t)) d\mu_0(x) ,
\end{equation}
where, now, the distribution of the diffusive species is described by a time-dependent measure $\mu(\cdot,t)$ with initial value $\mu_0$. Evidently, the regular form (\ref{dZ9cNP}) of the pushforward operation is recovered when $d\mu_0(x) = \rho_0(x) dx$ and $d\mu(y,t) = \rho(x,t) dy$ for absolutely continuous densities $\rho_0$ and $\rho$, respectively. In addition, expressing (\ref{yC2krR}) weakly, we obtain
\begin{equation}\label{VV8ZwT}
    \int_\Omega
        \partial_t \varphi(\varphi^{-1}(y,t),t)\cdot \xi(y)
    \, dy
    =
    \int_\Omega
        v(y,t) \cdot \xi(y)
    \, dy ,
\end{equation}
where $\xi$ are test functions {\black satisfying the boundary} condition
\begin{equation}\label{qS9jcJ}
    \xi \cdot n = 0 .
\end{equation}
Inserting the mobility law (\ref{rVre5S}) into (\ref{VV8ZwT}), gives
\begin{equation}
    \int_\Omega
        \partial_t \varphi(\varphi^{-1}(y,t),t) \cdot \xi(y)
    \, dy
    =
    \int_\Omega
        \big(
            \rho(y,t) u(y,t)
            -
            \kappa \nabla\rho(y,t)
        \big)
        \cdot \xi(y)
    \, dy .
\end{equation}
An integration by parts further gives
\begin{equation}
    \int_\Omega
        \partial_t \varphi(\varphi^{-1}(y,t),t) \cdot \xi(y)
    \, dy
    =
    \int_\Omega
        \big(
            u(y,t) \cdot \xi(y)
            +
            \kappa \, {\rm div} \, \xi(y)
        \big)
        \rho(y,t)
    \, dy ,
\end{equation}
where we have used (\ref{qS9jcJ}). As in the case of the pushforward operation, this weak reformulation of the mobility law admits the extension to measures
\begin{equation}\label{7M4w9Q}
    \int_\Omega
        \partial_t \varphi(\varphi^{-1}(y,t),t) \cdot \xi(y)
    \, dy
    =
    \int_\Omega
        \big(
            u(y,t) \cdot \xi(y)
            +
            \kappa \, {\rm div} \, \xi(y)
        \big)
    \, d\mu(y,t) ,
\end{equation}
suitable to approximation by particle methods.  Eqs.~(\ref{UjA9vT}) and (\ref{7M4w9Q}) may be thought as jointly defining an evolution for both the measure $\mu$ and the transport map $\varphi$. In particular, we note that the preceding extension to measures requires consideration of the transport map $\varphi$ as an additional unknown of the problem.

\subsection{Time discretization}
\label{sec:TimeDiscretization}

We begin by approximating problem (\ref{eq:TD:RV}) in time. To this end, let $t_0=0 < t_1 <$ $\cdots$ $<t_N = T$ be a discretization of the time interval {\black $(0,T)$}. The goal is {\black now} to determine corresponding discrete approximations $\rho_0$, $\rho_1$, $\dots$, $\rho_N$ of the densities and discrete approximations $\varphi_{0\to 1}$, $\varphi_{1\to 2}$, $\dots$, $\varphi_{N-1\to N}$ of the incremental transport maps.

We begin by noting that eq.~(\ref{eq:TD:Diff1}) may be regarded as an evolution equation in which the operator is the sum of two operators: advection and diffusion. This additive structure immediately suggests splitting the time integration into advective and diffusive fractional steps. The advective fractional step is governed by the advection equation
\begin{equation}
    \partial_t \rho + \nabla\cdot(\rho u) = 0,
\end{equation}
which follows formally from (\ref{eq:TD:Diff1}) by setting $\kappa = 0$. Conveniently, the advection equation can be solved exactly. Thus, let $\chi : \Omega \times [0,T] \to \Omega$ be the flow corresponding to the given advection velocity field $u$, i.~e.,
\begin{equation}
    \partial_t \chi(x,t) = u(\chi(x,t), t) .
\end{equation}
Then, the incremental advection of the density $\rho_k$ at time $t_k$ is the pushforward ($_\#$) of $\rho_k$ through $\chi_{k\to k+1}$
\begin{equation}
    \rho_{k+1} = {\chi_{k\to k+1}}_{\#} \rho_k ,
\end{equation}
where
\begin{equation}
    \chi_{k\to k+1}(y) = \chi(\chi^{-1}(y,t_k),t_{k+1})
\end{equation}
is the incremental advection mapping.

Next, we turn to the diffusive fractional step, governed by the diffusive equation
\begin{equation}
    \partial_t \rho = \kappa \Delta \rho,
\end{equation}
formally obtained from (\ref{eq:TD:Diff1}) by setting $u = 0$, and to the temporal discretization of the diffusion step. We characterize the incremental evolution of the density by means of the JKO functional
\begin{equation}\label{eq:TD:F}
    F(\varphi_{k\to k+1})
    =
    \frac{1}{2}
    \frac{d_W^2(\rho_k,\rho_{k+1})}{t_{k+1}-t_k}
    +
    \int_\Omega \kappa \rho_{k+1} \log\rho_{k+1} \, dx ,
\end{equation}
where $d_W$ is the Wasserstein distance between measures and
\begin{equation}
    \rho_{k+1} \circ \varphi_{k\to k+1}
    =
    \rho_k/\det\big(\nabla\varphi_{k\to k+1}\big) .
\end{equation}
The incremental transport map $\varphi_{k\to k+1}$ then follows from the minimum problem
\begin{equation}\label{eq:TD:Fmin}
    F(\varphi_{k\to k+1}) \to \min!
\end{equation}
We may formally verify that this formulation indeed supplies a time discretization of the mobility law (\ref{rVre5S}) in the limit of zero advection. Thus, taking variations we obtain (cf.~Appendix, eq.~(\ref{diffeWdist}), also \cite{Villani2003})
\begin{equation}\label{eq:TD:DFmin}
\begin{split}
    \langle DF(\varphi_{k\to k+1}), \xi_{k+1} \rangle
    & =
    \frac{1}{t_{k+1}-t_k}
    \int_\Omega \langle x-\varphi_{k+1\to k}(x), \xi_{k+1}(x) \rangle
    \rho_{k+1}(x) \, dx
    \\ & -
    \int_\Omega \kappa (\log\rho_{k+1} + 1) \nabla\cdot(\rho_{k+1}\xi_{k+1}) \, dx ,
    \\ & =
    \frac{1}{t_{k+1}-t_k}
    \int_\Omega \langle x-\varphi_{k+1\to k}(x), \xi_{k+1}(x) \rangle
    \rho_{k+1}(x) \, dx
    \\ & +
    \int_\Omega \kappa \nabla \rho_{k+1}\cdot \xi_{k+1} \, dx ,
\end{split}
\end{equation}
where we write
\begin{equation}
    \varphi_{k+1\to k}
    =
    \varphi_{k\to k+1}^{-1} ,
\end{equation}
and the variations satisfy the boundary condition
\begin{equation}
    \xi_{n+1} \cdot n = 0 .
\end{equation}
Enforcing stationarity for all admissible variations yields
\begin{equation}\label{eq:TD:RVD3}
    \rho_{k+1}(x)
    \frac{x-\varphi_{k+1\to k}(x)}{t_{k+1}-t_k}
    =
    -
    \kappa \nabla \rho_{k+1}(x) ,
\end{equation}
which is indeed a time discretization of (\ref{rVre5S}).

We note that the density $\rho_{k+1}$ still enters (\ref{eq:TD:DFmin}) nonlinearly and differentiated. In order to facilitate particle approximations, we perform an integration by parts leading to the weak stationarity condition
\begin{equation}\label{eq:TD:DFmin2}
    \int_\Omega
        \rho_{k+1}(y)
        \frac{y-\varphi_{k+1\to k}(y)}{t_{k+1}-t_k}
        \cdot
        \xi(y)
    \, dy
    =
    \int_\Omega
        \kappa \rho_{k+1}(y)
        \nabla \cdot \xi(y)
    \, dy .
\end{equation}
We verify that, indeed, the density $\rho_{k+1}$ now enters this weak stationarity condition linearly and undifferentiated.

Likewise, the incremental mass conservation relation (\ref{eq:TD:RV1}) must now be understood in a weak or distributional sense, i.~e., as the requirement that
\begin{equation}\label{eq:OT:CDsol3}
    \int \rho_k(x) \eta(x) \, dx
    =
    \int \rho_{k+1}(y)
    \eta\big(\varphi^{-1}_{k\to k+1}(y)\big) dy,
\end{equation}
for all test functions $\eta$. Again we verify that the density $\rho_{k+1}$ now enters this weak form of the transport equation linearly and undifferentiated.

\subsection{Spatial discretization}
\label{sec:SpaceDiscretization}

Next, we turn to the question of spatial discretization of the weak semi-discrete Fick's law (\ref{eq:TD:DFmin2}) and the weak form of the transport equation (\ref{eq:OT:CDsol3}). The structure of the latter reveals the need for two types of approximations: i) the discretization of the density $\rho_{k+1}$, and ii) the discretization of the incremental transport map $\varphi_{k+1\to k}$. We consider these two approximations in turn.

\subsubsection{Spatial discretization of the density.}

As already noted, the density $\rho_{k+1}$ enters (\ref{eq:TD:DFmin2}) and (\ref{eq:OT:CDsol3}) linearly and undiffrentiated. Therefore, a natural and computationally convenient spatial discretization may be effected by considering mass densities of the {\sl particle} type
\begin{equation}\label{eq:OT:Rhok}
    \rho_{h,k}(x)
    =
    \sum_{p=1}^M m_{p,k} \delta\big(x- x_{p,k}\big),
\end{equation}
where $x_{p,k}$ represents the position at time $t_k$ of a \emph{material-point} of mass $m_p${\black,} $\delta\big(x- x_{p,k}\big)$ is the Dirac-delta distribution centered at $x_{p,k}$, and $M$ is the number of material-points. For discrete mass distributions of the form (\ref{eq:OT:Rhok}), eq.~(\ref{eq:OT:CDsol3}) reduces to
\begin{equation}
    \sum_{p=1}^M m_{p,k} \eta( x_{p,k})
    =
    \sum_{p=1}^M m_{p,k+1} \eta( x_{p,k}),
\end{equation}
which must be satisfied for all test functions $\eta$. Hence,
\begin{equation}\label{eq:OT:CM}
    m_{p,k} = m_{p,k+1} = m_p,
\end{equation}
i.~e., the material-points must have constant mass, and representation (\ref{eq:OT:Rhok}) reduces to
\begin{equation}\label{eq:OT:Rhok2}
    \rho_{h,k}(x)
    =
    \sum_{p=1}^M m_p \delta\big(x- x_{p,k}\big) ,
\end{equation}
with constant $\{m_p,\ p = 1,\dots,M\}$. Thus, the weak and transport reformulation of the problem results trivially in {\sl exact mass conservation}, simply by keeping the mass of all particles constant.

\subsubsection{Spatial discretization of the incremental transport map.}

A full spatial discretization additionally requires the interpolation of the incremental transport map $\varphi_{k\to k+1}$. Since $\varphi_{k\to k+1}$ and its variations enter the governing equations (\ref{eq:TD:DFmin2}) and (\ref{eq:OT:CDsol3}) in differential form, the attendant interpolation must be conforming. To this end, we start by simplifying equation (\ref{eq:TD:DFmin2}) by recourse to the change of variables
\begin{subequations}
\begin{align}
    &
    x = \varphi_{k+1\to k}(y) ,
    \\ &
    y = \varphi_{k \to k+1}(x) ,
    \\ &
    \xi(x) = \eta(\varphi_{k \to k+1}(x)) ,
    \\ &
    \eta(y) = \xi(\varphi_{k+1\to k}(y)) ,
\end{align}
\end{subequations}
with the result
\begin{equation}\label{eq:SD:WR1}
    \int_\Omega
        \rho_k(x)
        \frac{\varphi_{k \to k+1}(x)-x}{t_{k+1}-t_k}
        \cdot
        \xi(x)
    \, dx
    =
    \int_\Omega
        \kappa \rho_k(x)
        {\rm tr}
        (
            \nabla\eta(\varphi_{k \to k+1}(x))
            \nabla\varphi_{k \to k+1}(x)
        )
    \, dx .
\end{equation}
We consider general linear interpolation schemes of the form
\begin{equation}\label{eq:SD:Varphih}
    \varphi_{h, k \to k+1}(x)
    =
    \sum_{a=1}^N x_{a,k+1} N_{a,k}(x) ,
\end{equation}
where $a$ indexes a {\sl nodal point set} within $\Omega$, $\{N_{a,k},\ a=1,\dots,N\}$ are {\black the} corresponding first-order consistent nodal shape functions at time $t_k${\black,} and $\{x_{a,k+1},\ a=1,\dots,N\} \equiv x_{k+1}$ is the array of unknown nodal coordinates at time $t_{k+1}$. Consistency here means, specifically, that the shape functions satisfy the identities
\begin{subequations}\label{eq:SD:Cons}
\begin{align}
    &
    \sum_{a=1}^N N_{a,k}(x) = 1, \label{eq:SD:Cons1}
    \\ &
    \sum_{a=1}^N x_{a,k} N_{a,k}(x) = x. \label{eq:SD:Cons2}
\end{align}
\end{subequations}
In addition, we interpolate the weight functions as
\begin{equation}
    \xi_h(x) = \sum_{a=1}^N \xi_a N_{a,k}(x)
\end{equation}
and
\begin{equation}
    \eta_h(y)
    =
    \sum_{a=1}^N \xi_a N_{a,k}(\varphi_{h,k+1\to k}(y)) .
\end{equation}
Inserting these representations into (\ref{eq:SD:WR1}) defines the fully-discrete problem
\begin{equation}\label{eq:SD:WR2}
\begin{split}
    &
    \int_\Omega
        \rho_k(x)
        \frac{\varphi_{h,k \to k+1}(x)-x}{t_{k+1}-t_k}
        \cdot
        \xi_h(x)
    \, dx
    = \\ &
    \int_\Omega
        \kappa \rho_k(x)
        {\rm tr}
        (
            \nabla\eta_h(\varphi_{h,k \to k+1}(x))
            \nabla\varphi_{h,k \to k+1}(x)
        )
    \, dx .
\end{split}
\end{equation}
However
\begin{equation}
    \nabla \eta_h(\varphi_{h,k \to k+1}(x))
    =
    \sum_{a=1}^N
    \xi_a \nabla N_{a,k}(x)
    \nabla \varphi_{h, k+1\to k}(\varphi_{h, k \to k+1}(x)) ,
\end{equation}
or
\begin{equation}
    \nabla \eta_h(\varphi_{h,k \to k+1}(x))
    =
    \sum_{a=1}^N
    \xi_a \nabla N_{a,k}(x)
    (\nabla\varphi_{h,k \to k+1})^{-1}(x) ,
\end{equation}
whereupon (\ref{eq:SD:WR2}) reduces to
\begin{equation}
\begin{split}
    &
    \int_\Omega
        \rho_k(x)
        \frac{1}{t_{k+1}-t_k}
        \left(
            \sum_{b=1}^N
            (x_{b,k+1}-x_{b,k}) N_{b,k}(x)
        \right)
        \cdot
        \left(
            \sum_{a=1}^N \xi_a N_{a,k}(x)
        \right)
    \, dx
    = \\ &
    \int_\Omega
        \kappa
        \rho_k(x)
        \left(
            \sum_{a=1}^N
            \xi_a \cdot \nabla N_{a,k}(x)
        \right)
    \, dx .
\end{split}
\end{equation}
Enforcing this equation for all nodal weights $\{\xi_a,\ a=1,\dots,N\}$ gives the system of equations
\begin{equation}\label{eq:SD:WR3}
    \int_\Omega
        \rho_k(x)
        \frac{1}{t_{k+1}-t_k}
        \left(
            \sum_{b=1}^N
            (x_{b,k+1}-x_{b,k}) N_{b,k}(x)
        \right)
        N_{a,k}(x)
    \, dx
    =
    \int_\Omega
        \kappa
        \rho_k(x)
        \nabla N_{a,k}(x)
    \, dx .
\end{equation}
Finally, inserting the discrete density (\ref{eq:OT:Rhok2}) into (\ref{eq:SD:WR3}) gives the particularly simple semi-discrete evolution equation
\begin{equation}\label{eq:SD:WR4}
    M_k \frac{x_{k+1} - x_k}{t_{k+1} - t_k}
    =
    f_k ,
\end{equation}
with {\black always positive definite} mass matrix
\begin{equation}
    M_{k,ab}
    =
    \sum_{p=1}^M
    m_p
    N_{a,k}(x_{p,k})
    N_{b,k}(x_{p,k}) ,
\end{equation}
abbreviated $M_k \equiv \{\{M_{k,ab}, \ a=1,\dots,N\},\ b=1,\dots,N\}$, and nodal flux array
\begin{equation}
    f_{k,a}
    =
    \sum_{p=1}^M
    m_p
    \kappa
    \nabla N_{a,k}(x_{p,k}) ,
\end{equation}
abbreviated $f_k \equiv\{f_{k,a},\ a=1,\dots,N\}$. We note that the mass matrices $M_k$ are symmetric.

\subsection{Update algorithm}
\label{sec:SpaceDiscretization}

Eq.~(\ref{eq:SD:WR4}) can be solved for the updated nodal coordinates, with the result
\begin{equation}\label{eq:SD:Update}
    x_{k+1}
    =
    x_k
    +
    (t_{k+1} - t_k) M_k^{-1} f_k .
\end{equation}
This scheme amounts to transporting the nodes ballistically over each time step with nodal velocities
\begin{equation}
    v_k
    =
    M_k^{-1} f_k .
\end{equation}
The resulting update for the diffusive fractional step is summarized in the following algorithm.

\begin{algorithm}
\begin{algorithmic}[1]
\STATE Initialization: Update nodal coordinates $x_{a,k}$, material-point coordinates $x_{p,k}$, by means of the advection map $\chi$.
\STATE Recompute shape functions $N_{a,k}(x_{p,k})$ and derivatives $\nabla N_{a,k}(x_{p,k})$ from advected nodal and material-point sets.
\STATE Compute mass matrix $M_k$ and nodal fluxes $f_k$:
\begin{equation}\nonumber
    M_{k,ab}
    =
    \sum_{p=1}^M
    m_p
    N_{a,k}(x_{p,k})
    N_{b,k}(x_{p,k})
\end{equation}
\begin{equation}\nonumber
    f_{k,a}
    =
    \sum_{p=1}^M
    m_p
    \kappa
    \nabla N_{a,k}(x_{p,k})
\end{equation}
\STATE Diffusive update nodal coordinates:
\begin{equation}\nonumber
    x_{k+1}
    =
    x_{k}
    +
    (t_{k+1} - t_k) M_k^{-1} f_k
\end{equation}
\STATE Update material-point coordinates:
\begin{equation}\nonumber
    x_{p,k+1} = \varphi_{h,k\to k+1}(x_{p,k})
\end{equation}
\STATE Update material-point volumes:
\begin{equation}\nonumber
    v_{p,k+1} = \det {\nabla \varphi_{h,k\to k+1}} \, v_{p,k}
\end{equation}
\STATE Update material-point mass densities:
\begin{equation}\nonumber
    \rho_{p,k+1} = \frac{m_p}{v_{p,k+1}}
\end{equation}
\STATE Reset $k\leftarrow k+1$. If $k=N$ exit. Otherwise go to (2).
\end{algorithmic}
\caption{{\sc Optimal transport scheme for the heat equation.}}
\label{Al:OTM:SF}
\end{algorithm}

This forward solution has the usual structure of explicit time-integration and updated-Lagrangian schemes. However, the scheme differs fundamentally from the classical schemes for advection-diffusion {\black equations}, which have a linear-space structure. {\black Thus}, the present scheme relies on the transport of nodes and mass points to evolve the system. In particular, all the finite kinematics of the motion, including advection and the mass density and volume updates, are {\sl geometrically exact}.

As {\black in} common practice with explicit schemes, the linear {\black solver} implied in (\ref{eq:SD:Update}) can be rendered trivial by {\sl mass lumping}. Thus, for instance, row-sum lumping yields the diagonal mass matrix
\begin{equation}
    M_{k,{\rm lumped}}
    =
    {\rm diag}
    \left \{
        \sum_{p=1}^M
        m_p
        N_{a,k}(x_{p,k}),
        \ a = 1,\dots,N
    \right \} .
\end{equation}
{\black Although this additional approximation strictly deviates from the minimum principle (\ref{eq:TD:F}) and, therefore, constitutes a {\sl variational crime}, it is nevertheless adopted for convenience in all numerical examples.}

We note that we are left with considerable latitude in the choice of shape functions. In the examples that follow, we choose to use meshfree max-ent shape functions \cite{ArroyoOrtiz2006} computed from the convected nodal coordinates. The continuous reconstruction of the shape functions has the effect of automatically reconnecting the material-points and the nodal set, with no need of remapping the information carried by the material-points. This property of the method is particularly convenient for nonlinear and history-dependent extensions of the method, in which the local material state often includes additional internal variable information.

From the computational point of view, a challenging aspect of methods based on optimal transportation theory is the periodic redefinition of the nodal and material point neighbor lists, which requires the use of efficient search algorithms \cite{LiHabbalOrtiz2010}. For large-scale problems, the reconstruction of neighbor-list adds overhead to the calculations. We mitigate this overhead by performing neighbor-list updates only when the distortion of a local neighborhood reaches a certain prespecified tolerance. The shape function update is performed at every step.

\section{Examples of Application}
\label{sec:3DExamples}

We present selected numerical simulations that showcase the range and robustness of the approach in examples involving pure diffusion, pure advection, and combined advection-diffusion. In all calculations, the initial nodal-point set is obtained through an auxiliary Delaunay triangulation of the support of the initial density, with material points then placed at the barycenters of the corresponding tetrahedra.

{\black The critical step for the advection step is determined by a Courant condition bases on particle velocity and mean-free path between particles. For the small particle velocities in the examples under consideration, diffusion controls the critical time step. The time-step requirement is estimated with the following relation
\begin{equation}\label{eq:stableTimeStep}
    \Delta t \ll \frac{\Delta x^2}{\kappa} \, .
\end{equation}
where $\Delta x$ is the minimum distance between two nodes belonging to the same material-point neighborhood. It is also possible to advance the advection and diffusive steps with different time steps using subcycling, but such extensions of the method are beyond the scope of the paper.}
The time step is recomputed frequently, in particular after the redefinition of the connectivity table. In all examples, diffusion and possibly advection take place concurrently within confined volumes. Neumann boundary conditions are enforced by controlling the position of the nodes. Thus, nodes exiting the domain are repositioned on the boundary. No constraints are applied to the material-points.

\subsection{Pure diffusion in a spherical volume}
\label{ssec:sphere}

\begin{table}[!h]
\begin{center}
\begin{tabular}{lccc}
  \hline
  Mesh    & Nodes & material-points & $h_{\rm min}$\\
  \hline
  Coarse  & 235    & 593      & 0.26 \\
  Medium  & 442    & 1,259    & 0.17 \\
  Fine    & 1,510  & 4,997    & 0.12 \\
  \hline
\end{tabular}
\caption{\footnotesize Meshes to be used for the convergence analysis.}
\label{table:meshSphere}
\end{center}
\end{table}
We consider mass diffusion in a spherical volume from an initial configuration in the form of a sphere of radius $R_0=1$ with uniform density $\rho_0=1$. The mass is discretized into particle sets of three different sizes ranging from $593$ to $4,997$ material points, cf.~Table~\ref{table:meshSphere}. Material-points are subsequently free to diffuse within a sphere of radius $R_1 = 7$.
\begin{figure}[htp!]
\begin{center}
    \subfigure[Initial configuration]
    {\includegraphics[width=0.31\textwidth]{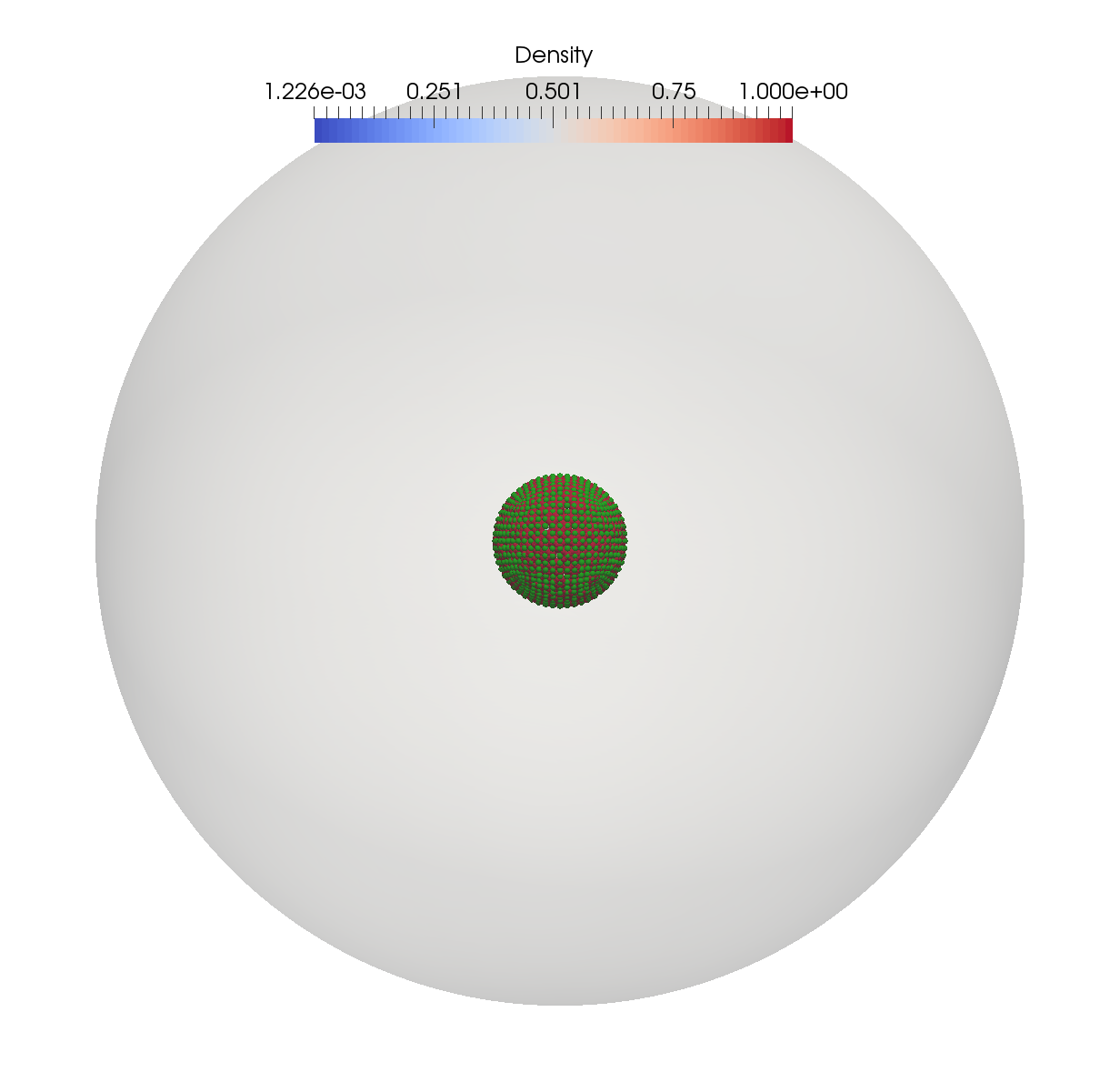}
    \label{fig:mpInitial}}
    \hskip 0.01\textwidth	
    \subfigure[Final configuration of material-points]
    {\includegraphics[width=0.31\textwidth]{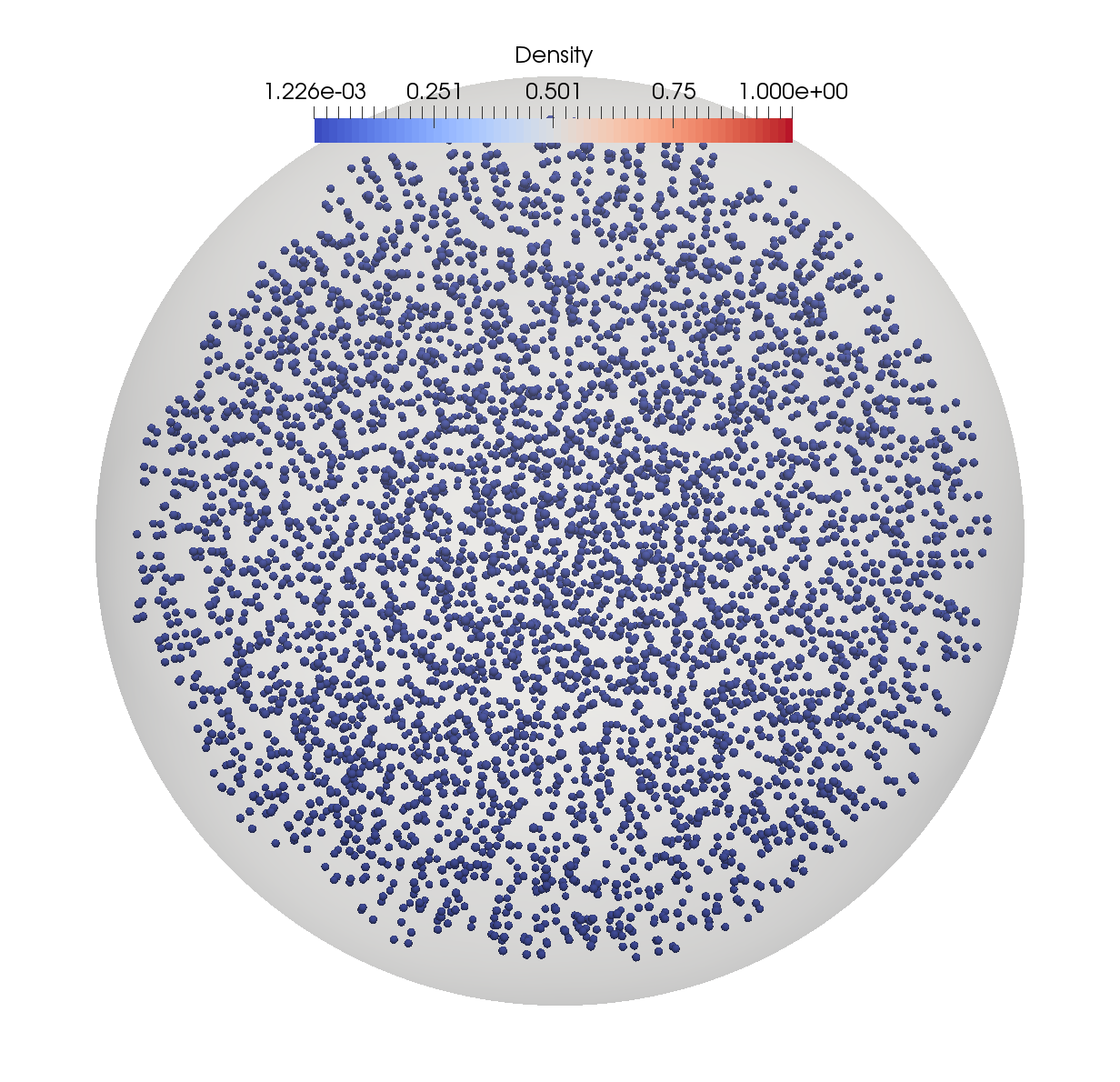}
    \label{fig:mpFinal}}
    \hskip 0.01\textwidth	
    \subfigure[Final configuration of nodes]
    {\includegraphics[width=0.31\textwidth]{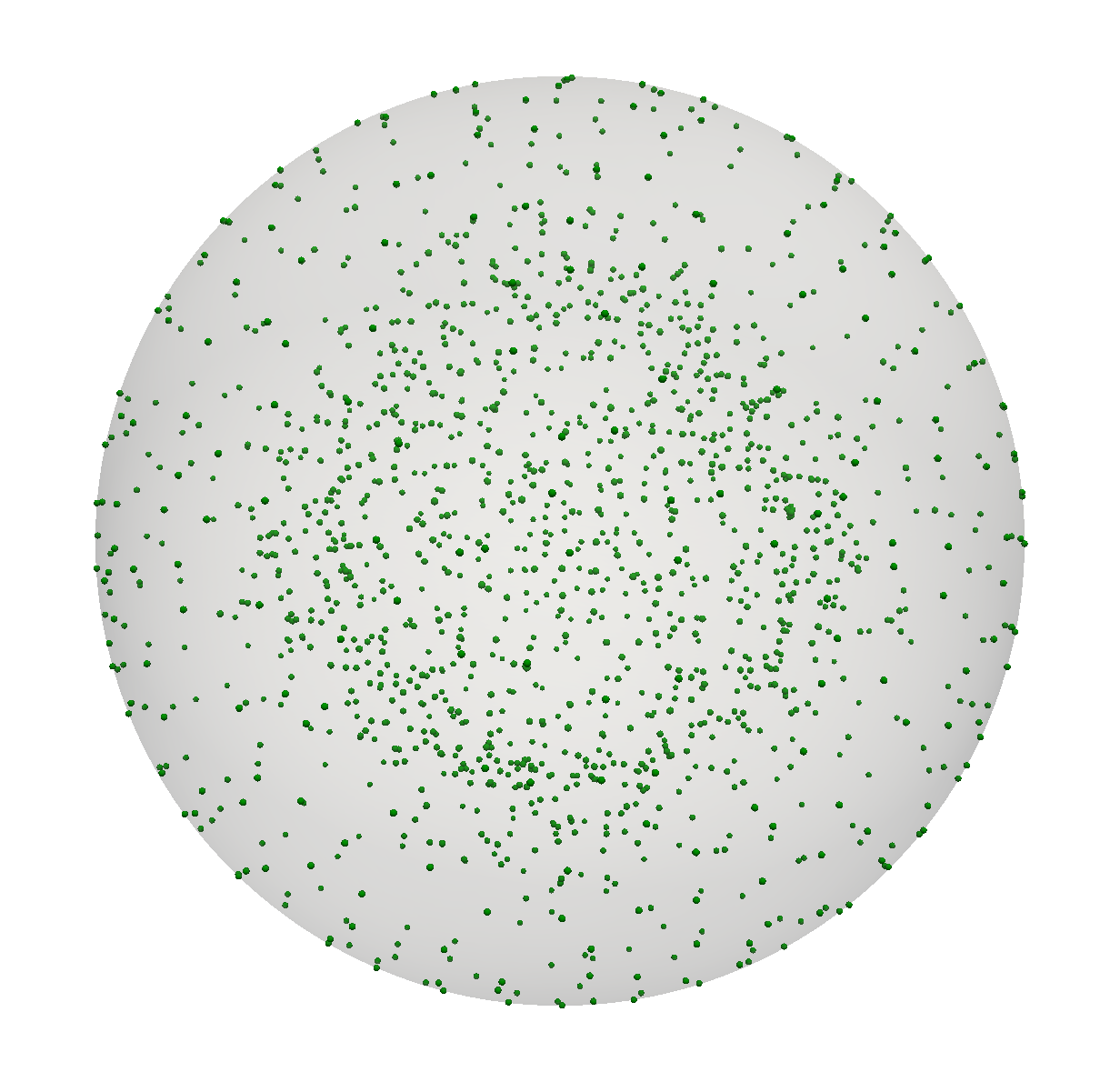}
    \label{fig:nodesFinal}}
	\caption{\footnotesize Expanding sphere problem, fine mesh. Configurations of the particle system, expanding from an initial spherical volume of radius $R_0=1$ within a spherical container of radius $R_f=7$. Each particle is colored according to its density, i.~e., the mass of the particle divided by its volume. Nodes are uniformly colored.}
	\label{fig:Sphere}
\end{center}
\end{figure}
\begin{figure}[htp!]
\begin{center}
    \subfigure[Density and dimensionless volume]
    {\includegraphics[width=0.45\textwidth]{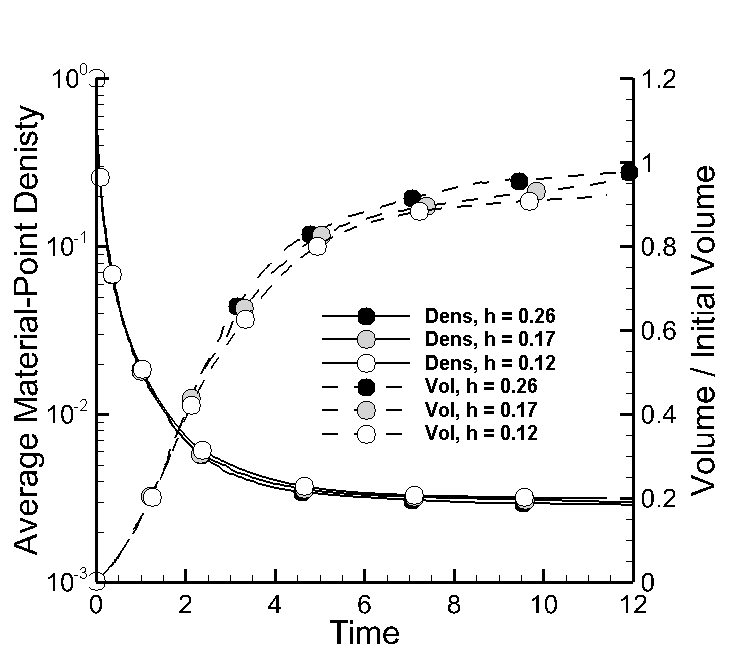}
    \label{fig:densityVolume}}
    \hskip 0.2cm	
    \subfigure[Radius]
    {\includegraphics[width=0.45\textwidth]{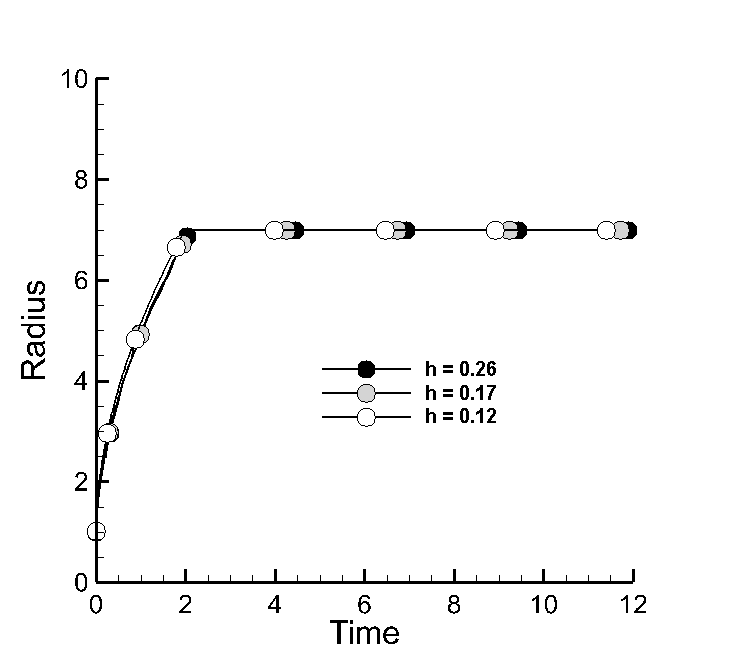}
    \label{fig:SphereRadius}}
	\caption{\footnotesize Expanding sphere problem. Time history of
    global variables for three different {\black discretizations}, cf. Table~\ref{table:meshSphere}. (a) Average density (in logarithmic scale) and dimensionless volume occupied by the material-points. (b) Maximum distance of the nodes from the center of the sphere.} \label{fig:SphereDensityRadius}
\end{center}
\end{figure}
For the fine mesh, Fig.~\ref{fig:mpInitial} shows the initial configuration of material-points (red) and nodes (green), whereas Figs.~\ref{fig:mpFinal}-\ref{fig:nodesFinal} show the final configuration attained by the material points and nodes, respectively. We note that the mass distribution evolves towards a uniform density in which the material points are uniformly distributed and carry ostensibly the same volume, Fig.~\ref{fig:mpFinal}. For the three particle-set sizes, Fig.~\ref{fig:densityVolume} displays the time history of the average density of the system and the volume occupied by the material-points. The expected monotonic trend towards uniform density is also clear from these figures. In terms of the minimum principle (\ref{eq:TD:F}), this trend may be understood as the result of a competition between entropy, which favors a uniform distribution, and the Wasserstein distance, which penalizes the movement of the particles. Fig.~\ref{fig:SphereRadius} shows the time history of the radius of the support of the particles, which experiences a steady expansion from its initial value until such time as the particles collide with the boundary of the container. Note that freezing the nodal points at the boundary of the sphere enforces condition \eqref{eq:TD:Diff2}. These results are nearly independent of the number of particles over the range considered, which demonstrates the fast convergence of the method and the high accuracy of relatively coarse approximations.

\subsection{Advection-diffusion in a hollow cylinder}
\label{ssec:hollowCylinder}

\begin{figure}[htp!]
\begin{center}
	\subfigure[$t=0.0$]
    {\includegraphics[width=0.45\textwidth]{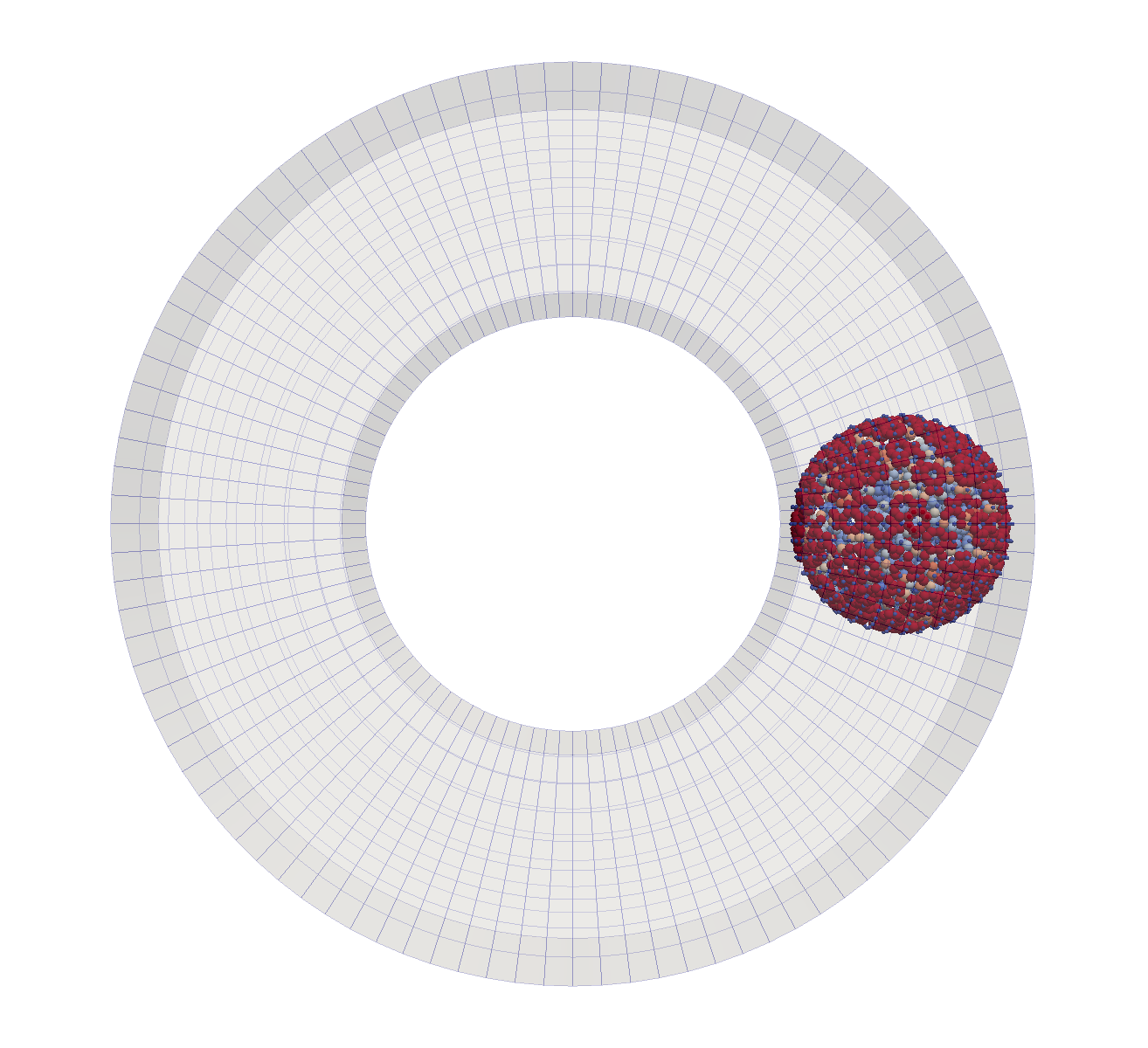}
    \label{fig:a1}}
	\subfigure[$t=0.3$]
    {\includegraphics[width=0.45\textwidth]{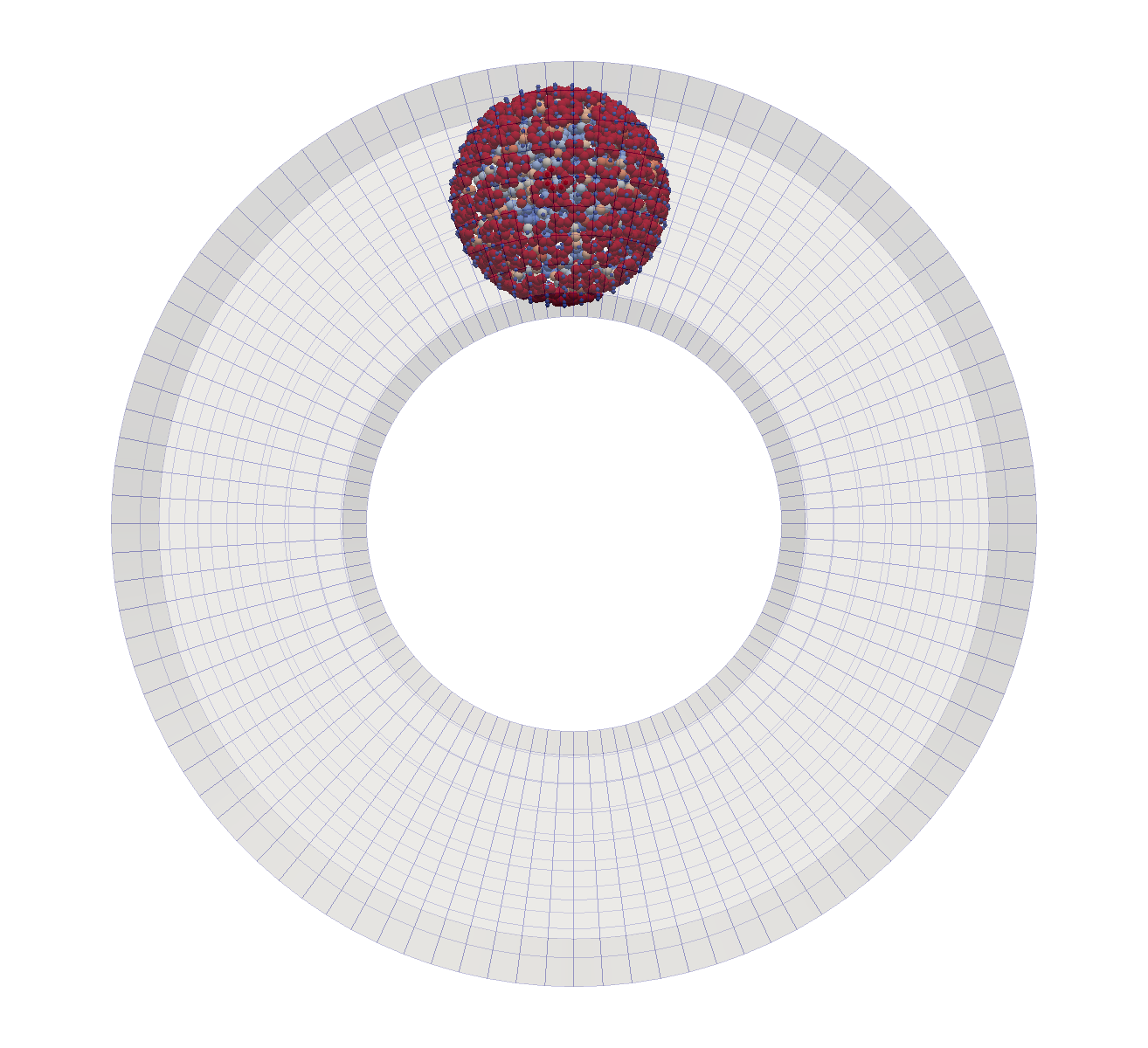}
    \label{fig:a2}}
	\subfigure[$t=3.0$]
    {\includegraphics[width=0.45\textwidth]{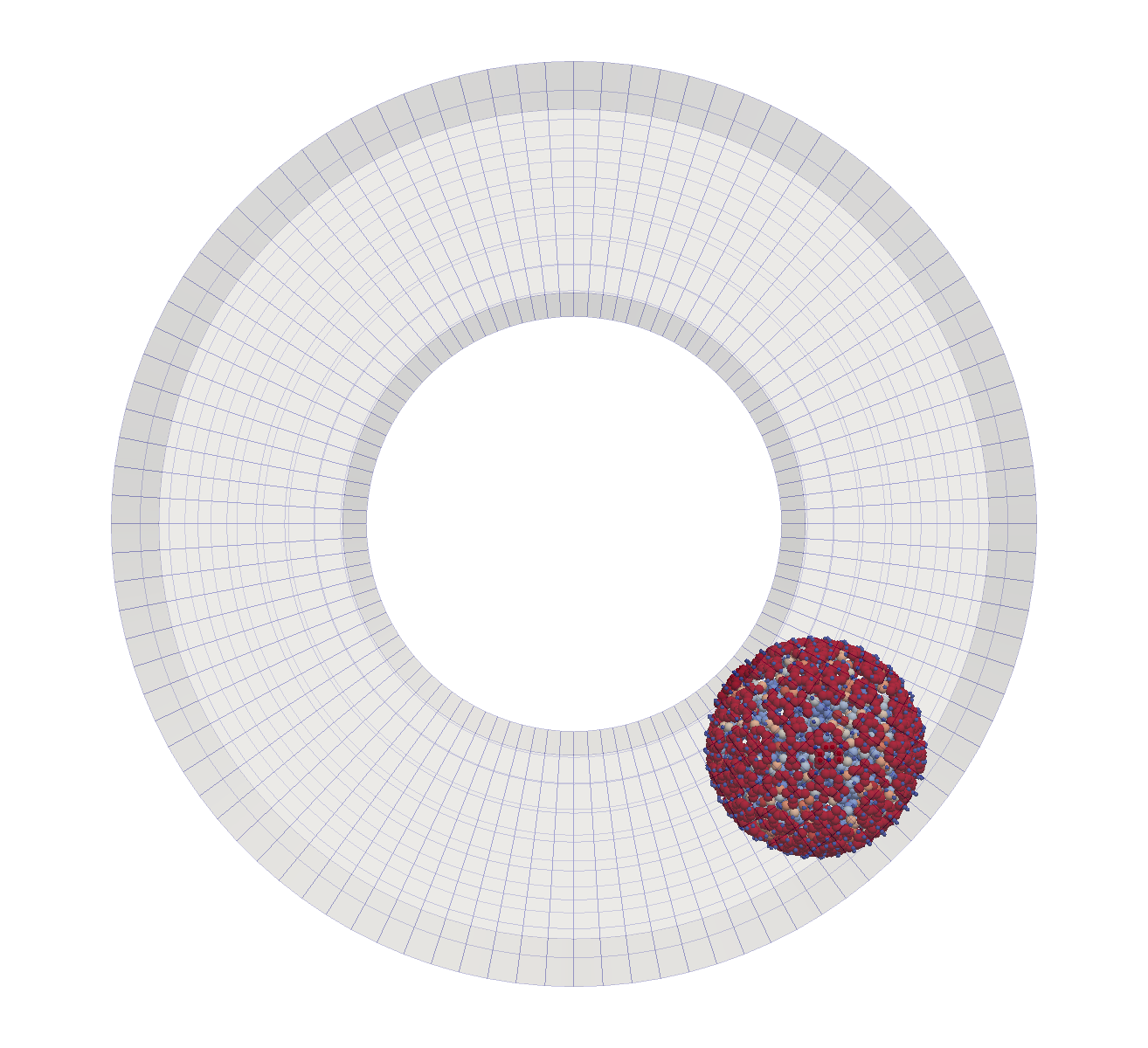}
    \label{fig:a3}}
	\subfigure[$t=3.0$, detail]
    {\includegraphics[width=0.45\textwidth]{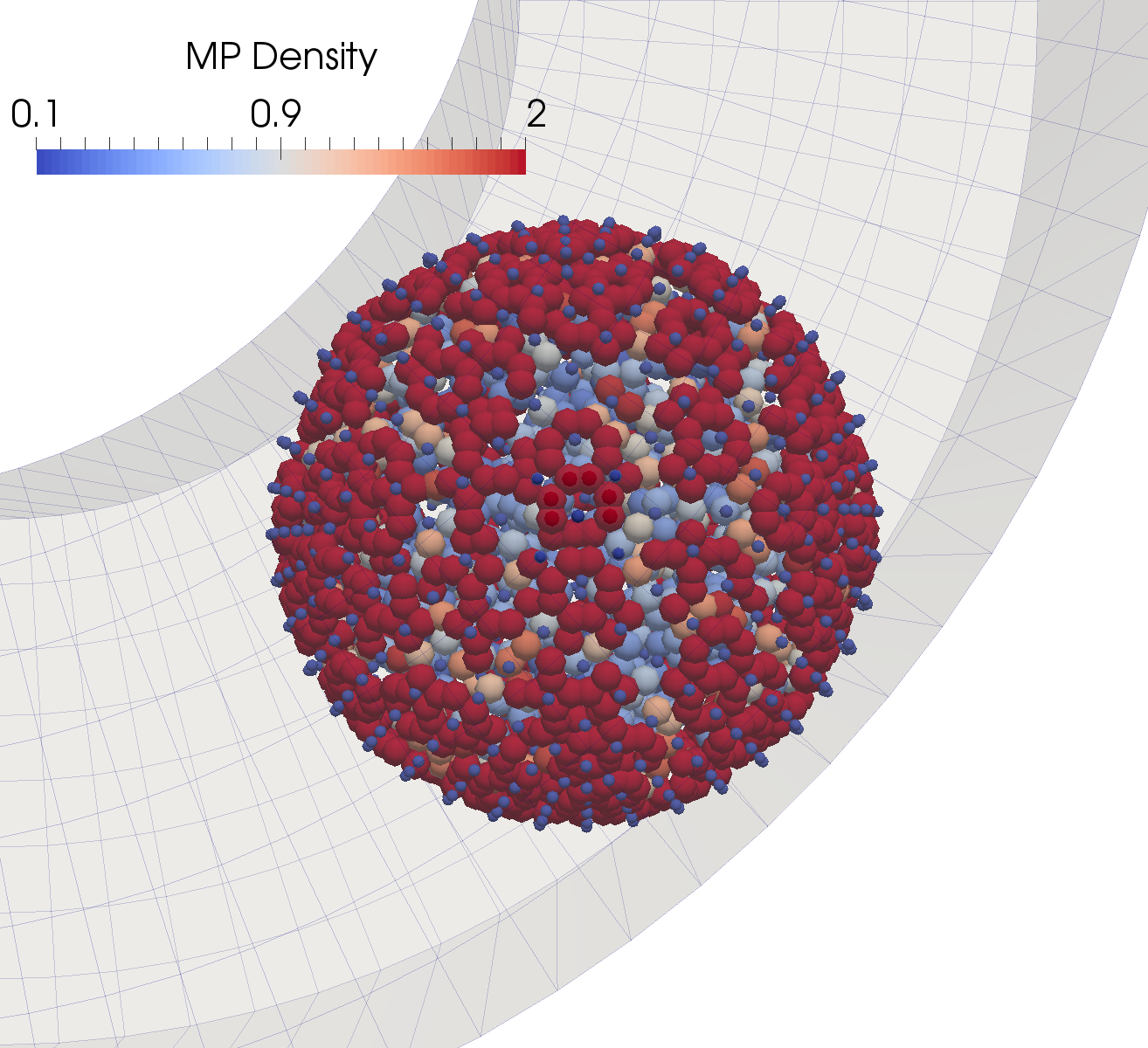}
    \label{fig:a4}}
	\caption{\footnotesize
    Pure advection in a circular channel of square cross section. The initial density is uniform over a spherical region spanning the cross section of the channel. The advection velocity field corresponds to a rigid-body motion at constant angular velocity. Snapshots of the material-point set at three successive times. The last image shows a close-up of the last configuration. Material points are colored according to the density, and the nodes are in blue.}
	\label{fig:CylinderAdvection}
\end{center}
\end{figure}

As a second example, we consider the problem of advection-diffusion in a circular {\black annulus with square cross section and periodic boundary conditions}, Fig.~\ref{fig:CylinderAdvection}a. The external radius of the {\black channel} is $0.5$ and the internal radius $0.25$. The initial mass density is uniform within a spherical region spanning the cross section of the channel. The density is discretized into 723 nodes and 2,215 material-points. The advection velocity field corresponds to a rigid-body motion at a constant angular velocity of 4. We consider the case of pure advection, $\kappa=0$, and the case of {\black combined advection-diffusion with diffusivity} $\kappa=0.001$.

Successive snapshots of the density evolution in the pure-advection problem spanning the short-term part of the solution are shown in Fig.~\ref{fig:CylinderAdvection}. Remarkably, the density function is advected {\sl exactly} by the algorithm with no numerical diffusion or spurious noise of any type. This type {\black of} behavior is, of course, a manifestation of the {\sl geometrically exact} character of the algorithm and is in sharp contrast to linear-space schemes, which experience difficulty in dealing with advection.

\begin{figure}[htp!]
\begin{center}
	\subfigure[$t=0.3$]
    {\includegraphics[width=0.45\textwidth]{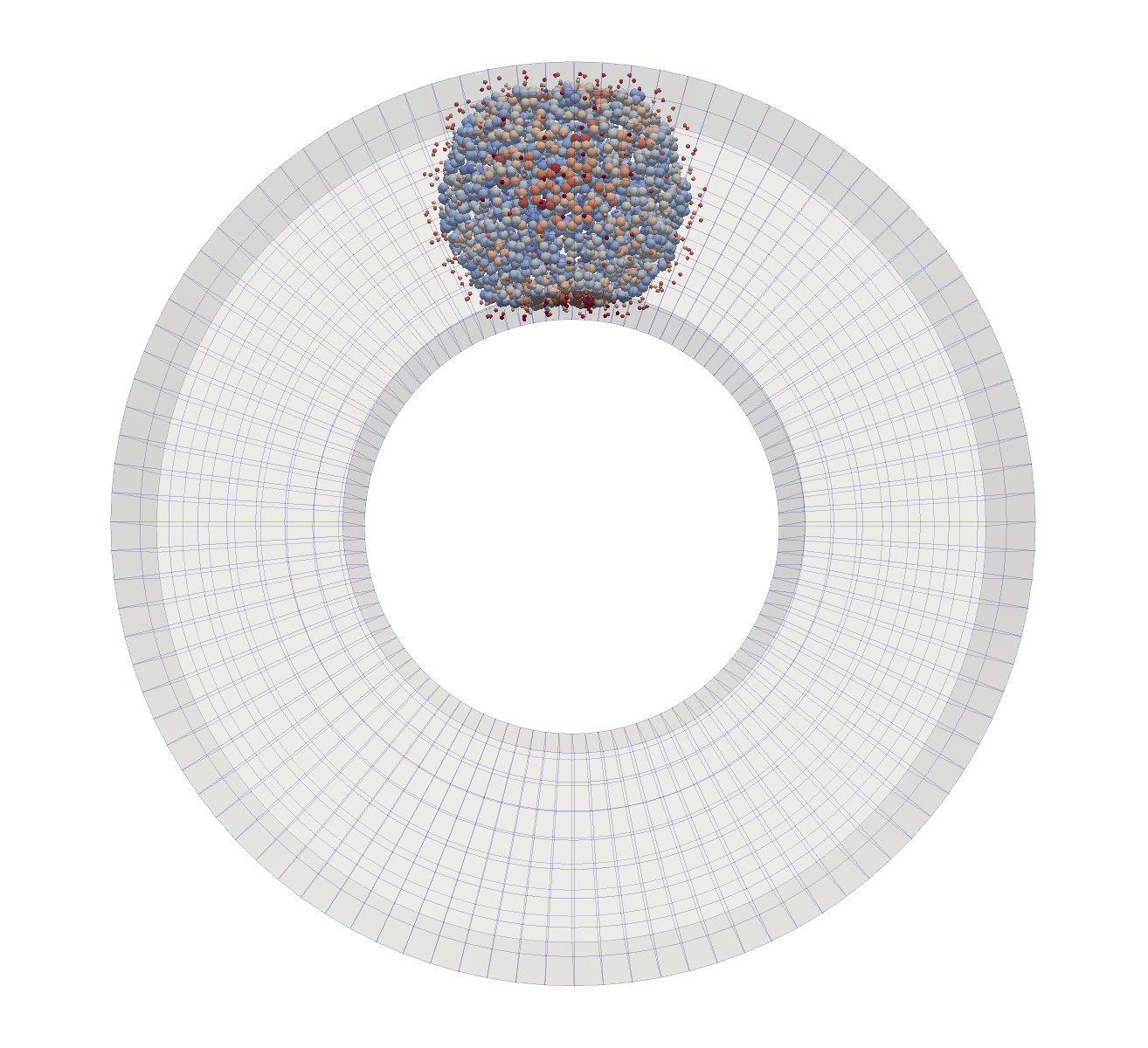}
    \label{fig:ad1}}
	\subfigure[$t=0.6$]
    {\includegraphics[width=0.45\textwidth]{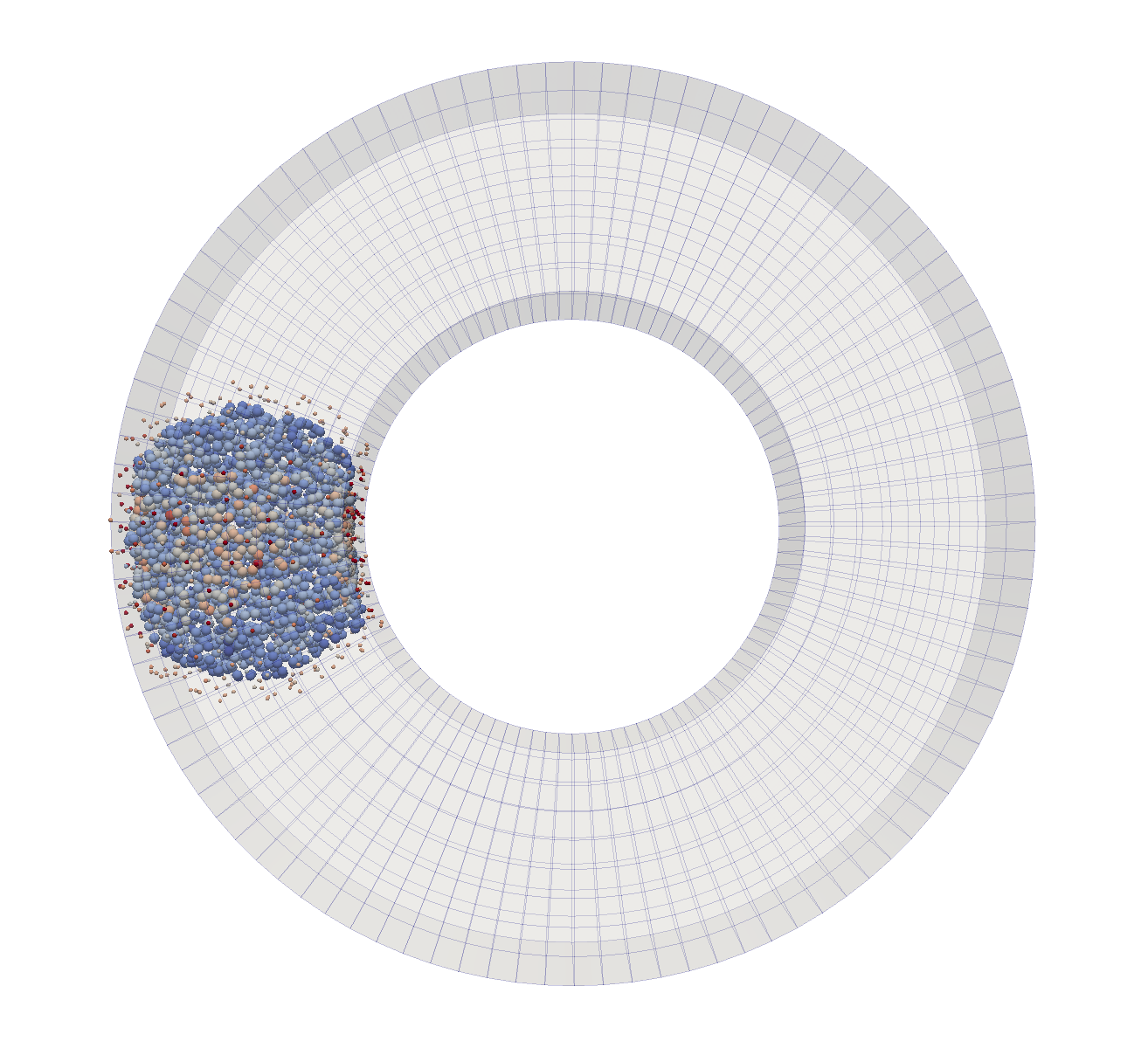}
    \label{fig:ad2}}
	\subfigure[$t=3.0$]
    {\includegraphics[width=0.45\textwidth]{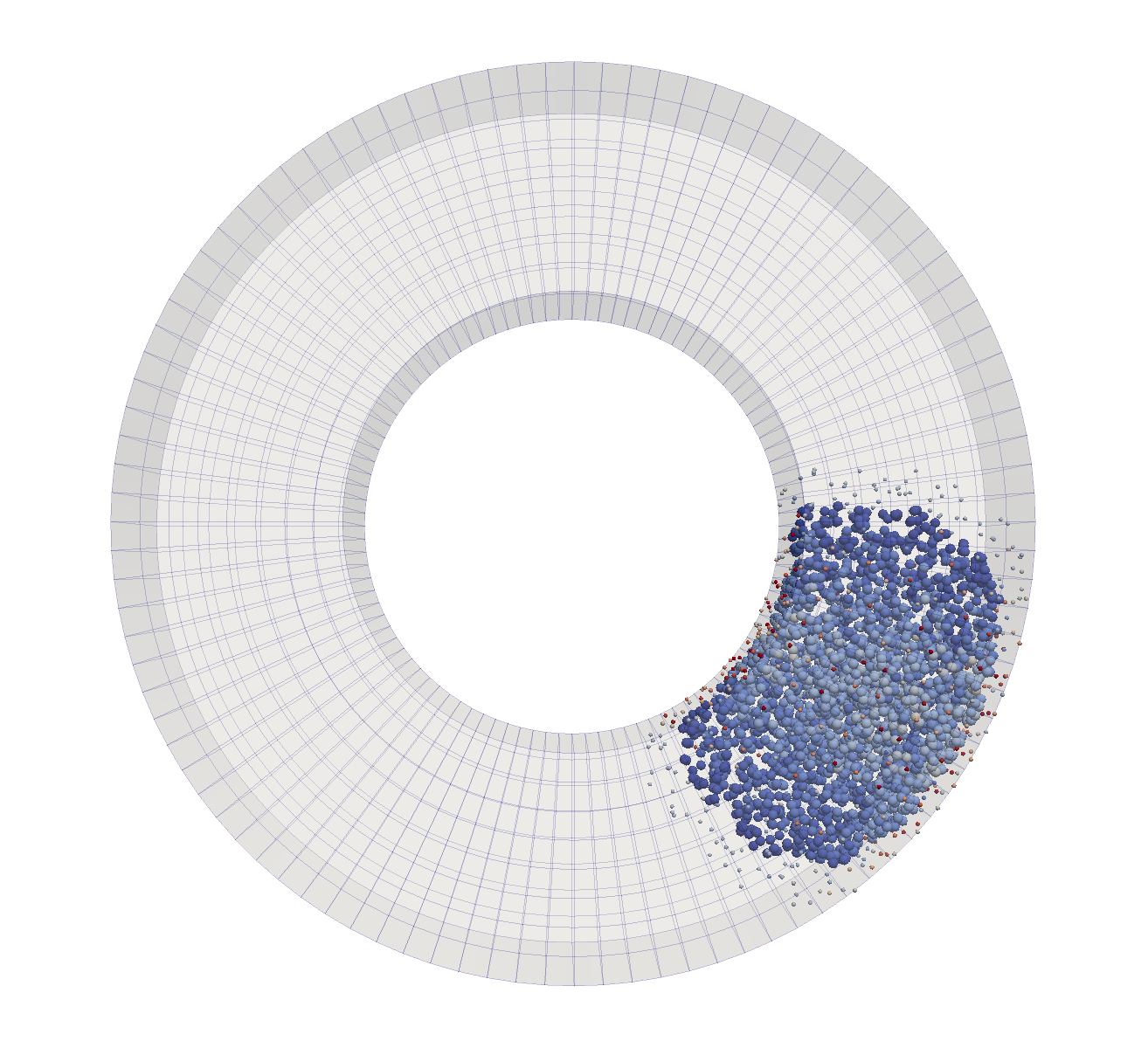}
    \label{fig:ad3}}
	\subfigure[$t=3.0$, detail]
    {\includegraphics[width=0.45\textwidth]{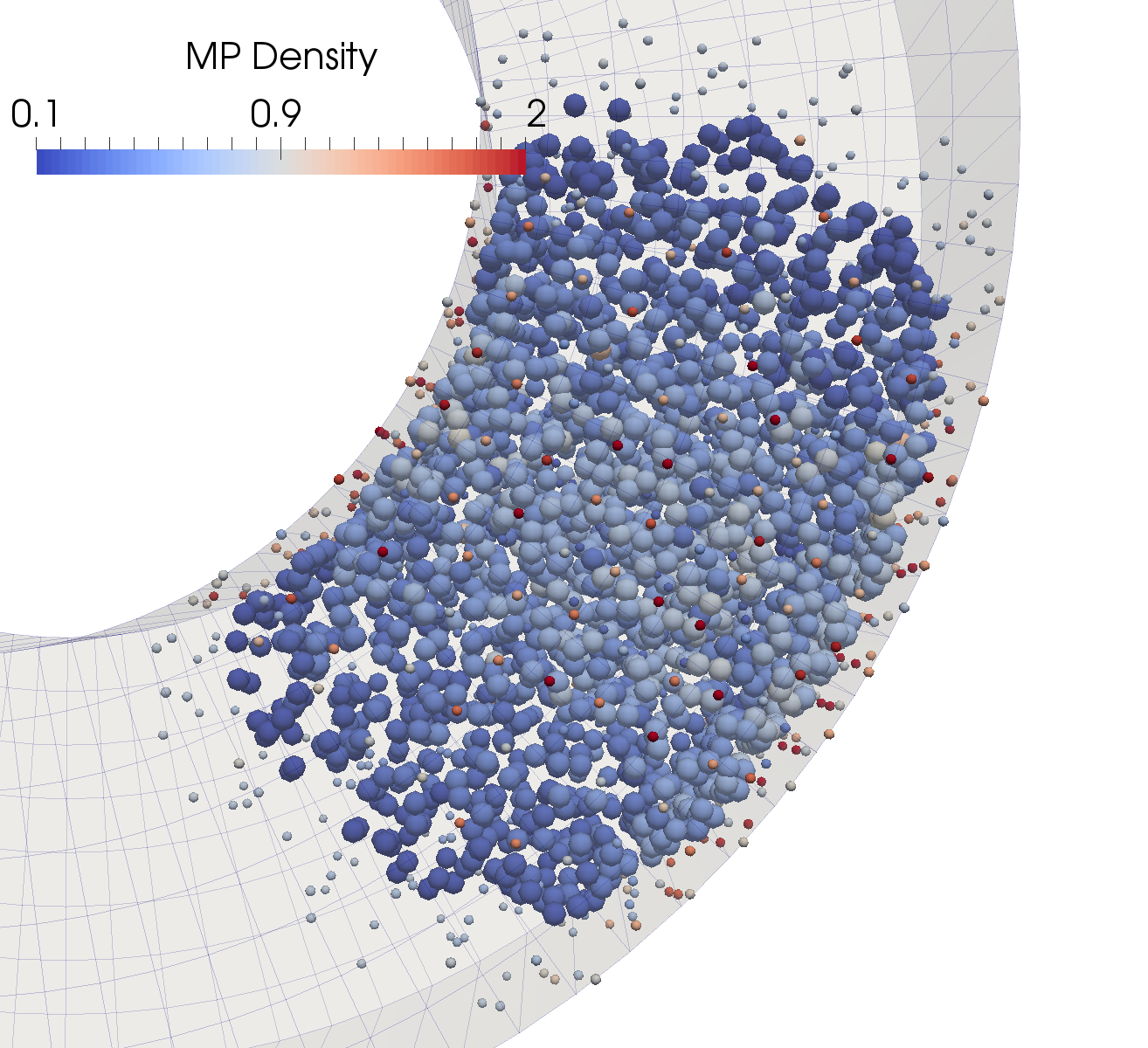}
    \label{fig:ad3}}
	\caption{\footnotesize
    Advection-diffusion in a circular channel of square cross section. The initial density is uniform over a spherical region spanning the cross section of the channel. The advection velocity field corresponds to a rigid-body motion at constant angular velocity. Snapshots of the material-point set at three successive times. The last image shows a close-up of the last configuration. Material points are colored according to the density, and the nodes are in blue.}
	\label{fig:CylinderAdvectionDiffusion}
\end{center}
\end{figure}

Fig.~\ref{fig:CylinderAdvectionDiffusion} collects successive snapshots of the advection-diffusion solution. In this case, a certain amount of diffusion tending to spread out the mass distribution is superimposed on the overall circular motion imparted by the advective flow. Eventually, the particles fill the entire channel and the mass density attains a uniform value. The robustness and stability with which the algorithm transports the mass particles and accounts for complex particle-particle and particle-boundary interactions is noteworthy.

\begin{figure}[htp!]
\begin{center}
	\includegraphics[width=0.45\textwidth]
    {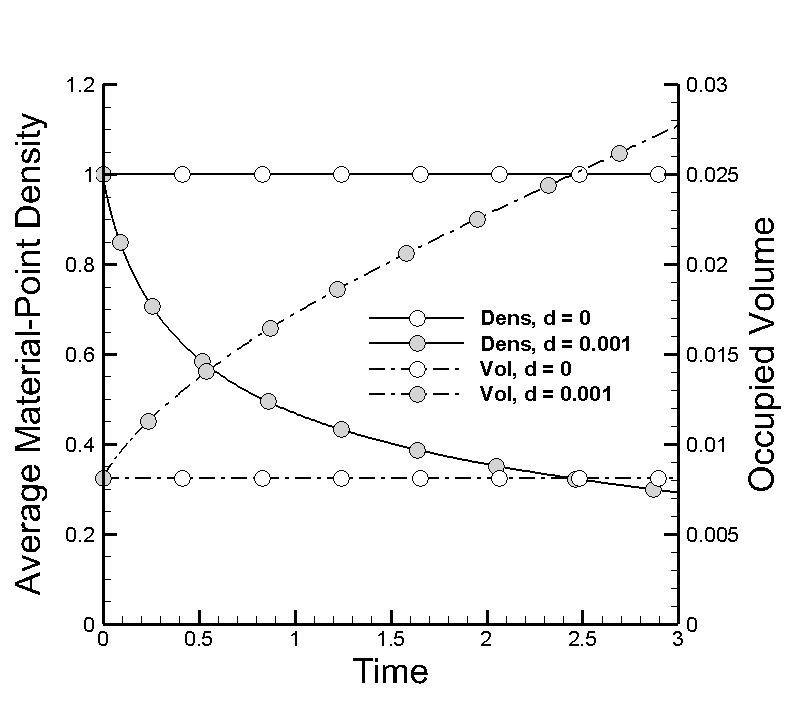}
	\caption{\footnotesize Advection-diffusion in a circular channel of square cross section. The initial density is uniform over a spherical region spanning the cross section of the channel. The advection velocity field corresponds to a rigid-body motion at constant angular velocity. Time evolution of the average material-point density and the volume occupied by the material-point set.}
	\label{fig:AdvectionDensityVolume}
\end{center}
\end{figure}

Fig.~\ref{fig:AdvectionDensityVolume} shows the time history of the average material-point density and of the volume occupied by the material-point set for the pure advection and the advection-diffusion problems. As expected, the average density {\black and} occupied volume remain exactly constant {\black for} pure-advection, again a manifestation of the geometrically exact nature of the method. In the advection-diffusion problem, {\black the} average density steadily decreases and, correspondingly, the occupied volume {\black increases with} time. The smoothness of the average quantities in time is noteworthy.

\begin{figure}[htp!]
\begin{center}
	\subfigure[$t=0.3$]
    {\includegraphics[width=0.45\textwidth]{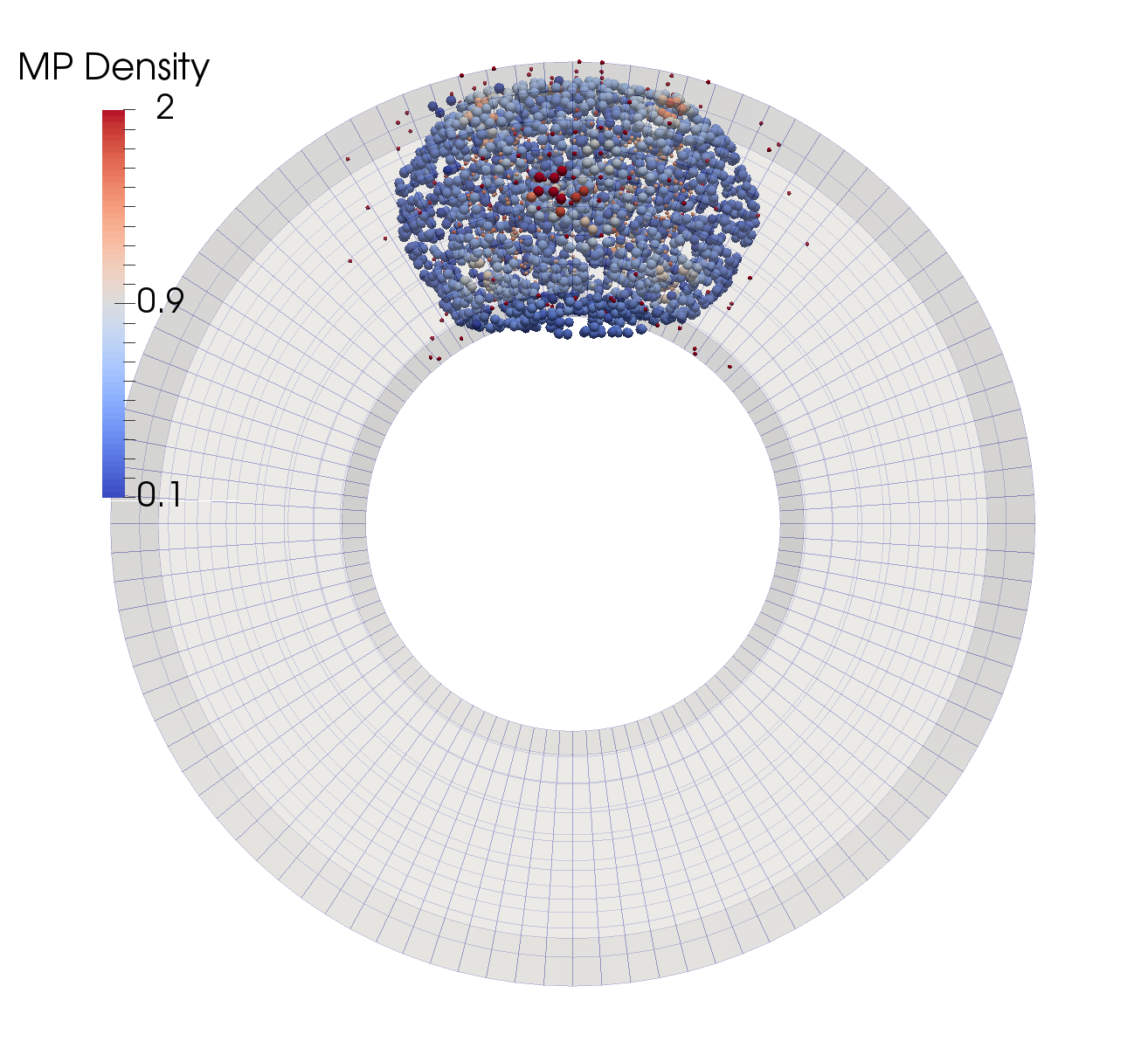}
    \label{fig:CylinderDiffusion1}}
	\subfigure[$t=0.6$, detail]
    {\includegraphics[width=0.45\textwidth]{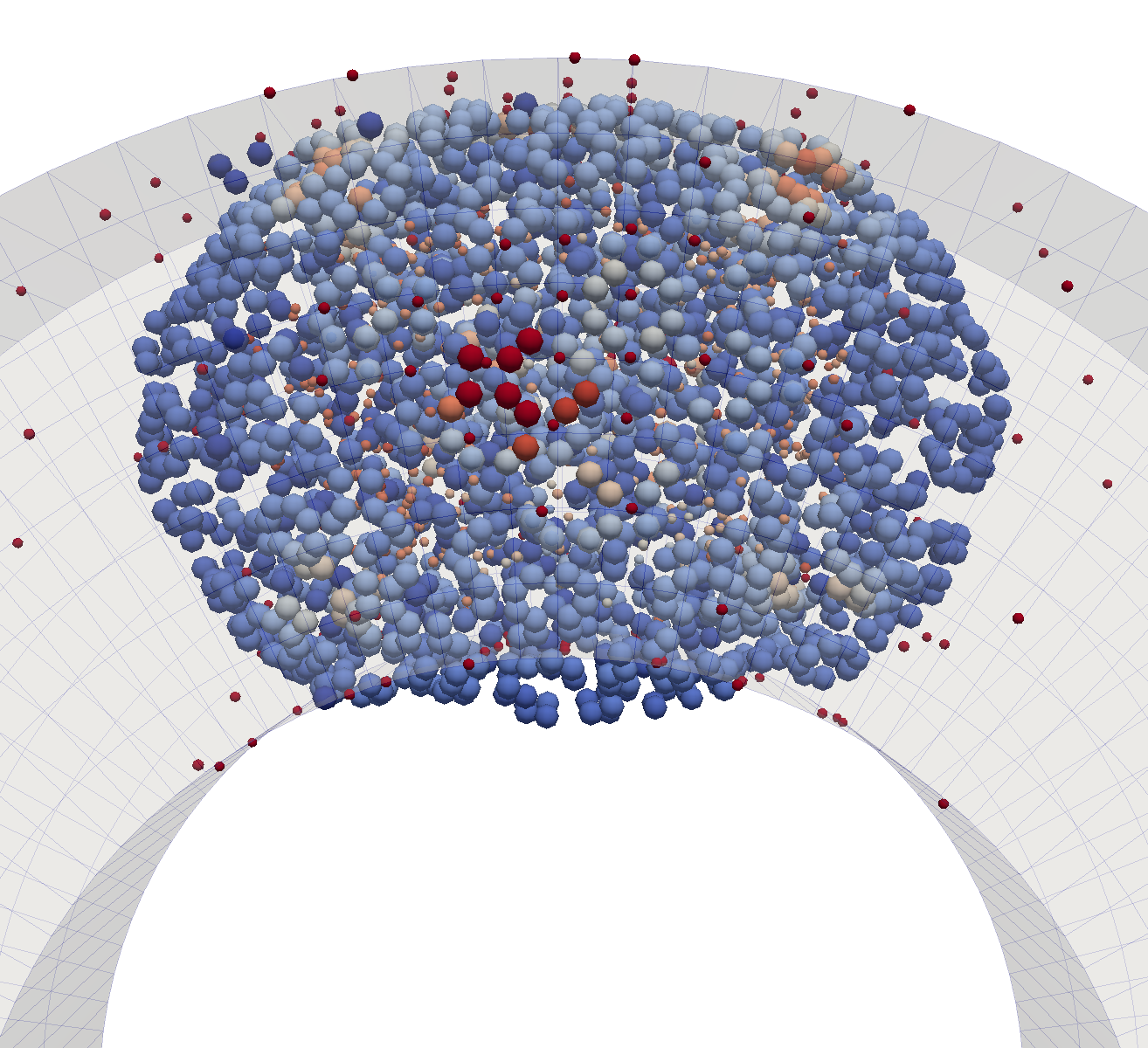}
    \label{fig:CylinderDiffusion2}}
  	\caption{\footnotesize
  Under-resolved advection-diffusion in a circular channel of square cross section exhibiting spill over the the material and nodal sets across the boundary.}
	\label{fig:CylinderDiffusion}
\end{center}
\end{figure}

A number of considerations of spatial and time resolution need to be carefully born in mind. Firstly, max-ent meshless shape functions such as formulated in \cite{ArroyoOrtiz2006} are defined through a convex programming problem that is sure to yield solutions when the domain is convex. If, contrariwise, the boundary has concave regions, the nodal point set must resolve the minimum radius of curvature $R$, or {\sl feature size}, for the max-ent shape functions to be well defined. Thus, we must have
\begin{equation}\label{eq:stableMeshSize}
    \Delta x \ll R \, ,
\end{equation}
which places spatial resolution requirements in terms of the geometry of the domain. In addition, as already noted the explicit nature of present formulation places restrictions on the time step, cf.~eq.~\ref{eq:stableTimeStep}. Combining (\ref{eq:stableTimeStep}) and (\ref{eq:stableMeshSize}) we additionally have
\begin{equation}\label{eq:stableTimeStep2}
    \Delta t \ll \frac{R^2}{\kappa} \, ,
\end{equation}
which places temporal resolution requirements in terms of the geometry of the domain. The consequences of insufficient resolution are illustrated in Fig.~\ref{fig:CylinderDiffusion}, corresponding to the problem of advection diffusion in a circular channel of square cross-section discussed in the foregoing with a larger diffusivity of $\kappa=0.01$, which tightens the resolution requirements. Thus, whereas the material and nodal point sets provide sufficient resolution for a diffusivity $\kappa=0.001$, they fail to do so for the larger diffusivity{\black ,} with the result that the material and nodal point sets spill over the concave part of the boundary and exit the domain of analysis.

\subsection{Advection and diffusion in a bucket}
\label{ssec:bucket}

\begin{figure}[htp!]
\begin{center}
	\subfigure[0.0]
    {\includegraphics[width=0.3\textwidth]{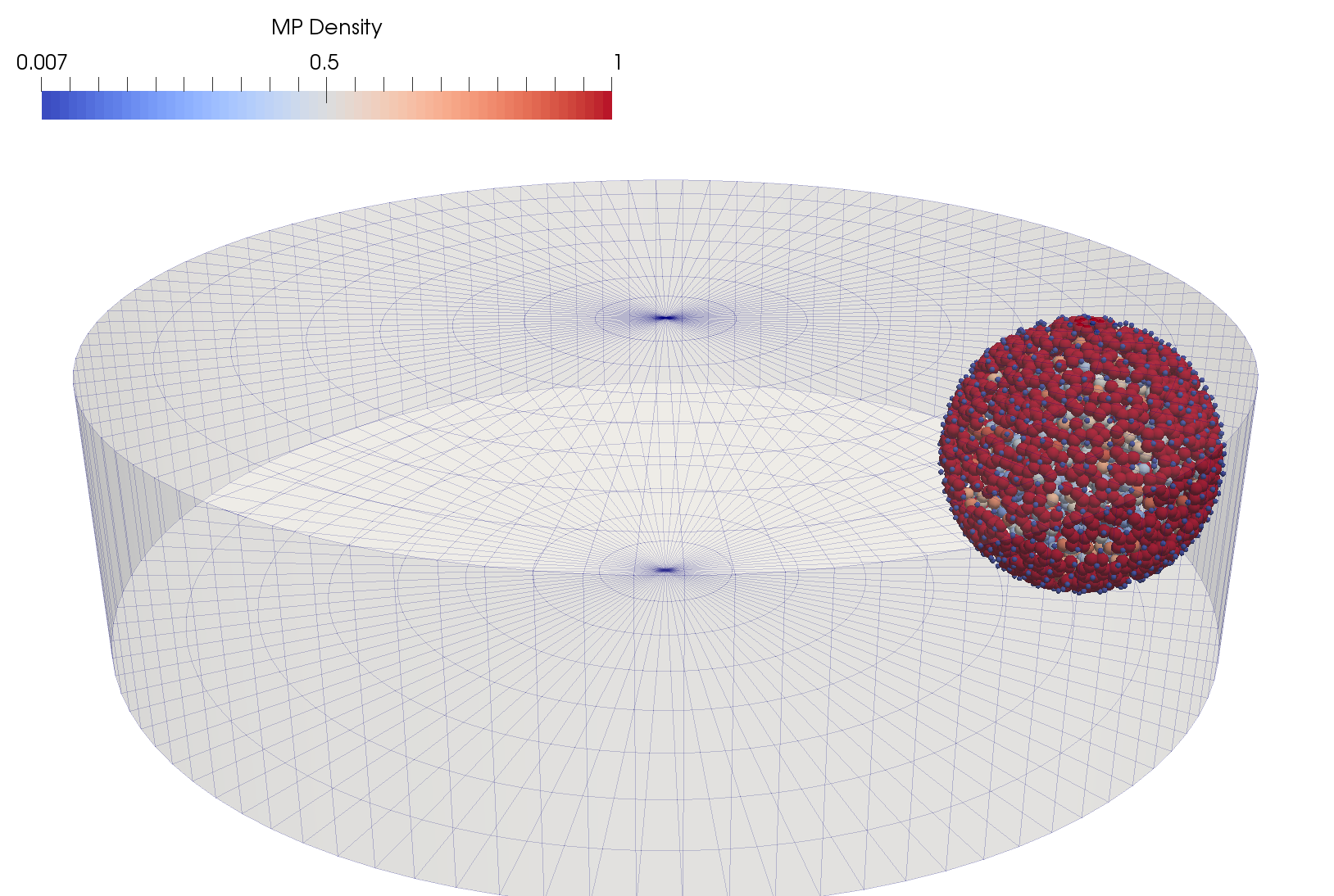}
    \label{fig:b1}}
	\subfigure[$t=0.1$]
    {\includegraphics[width=0.3\textwidth]{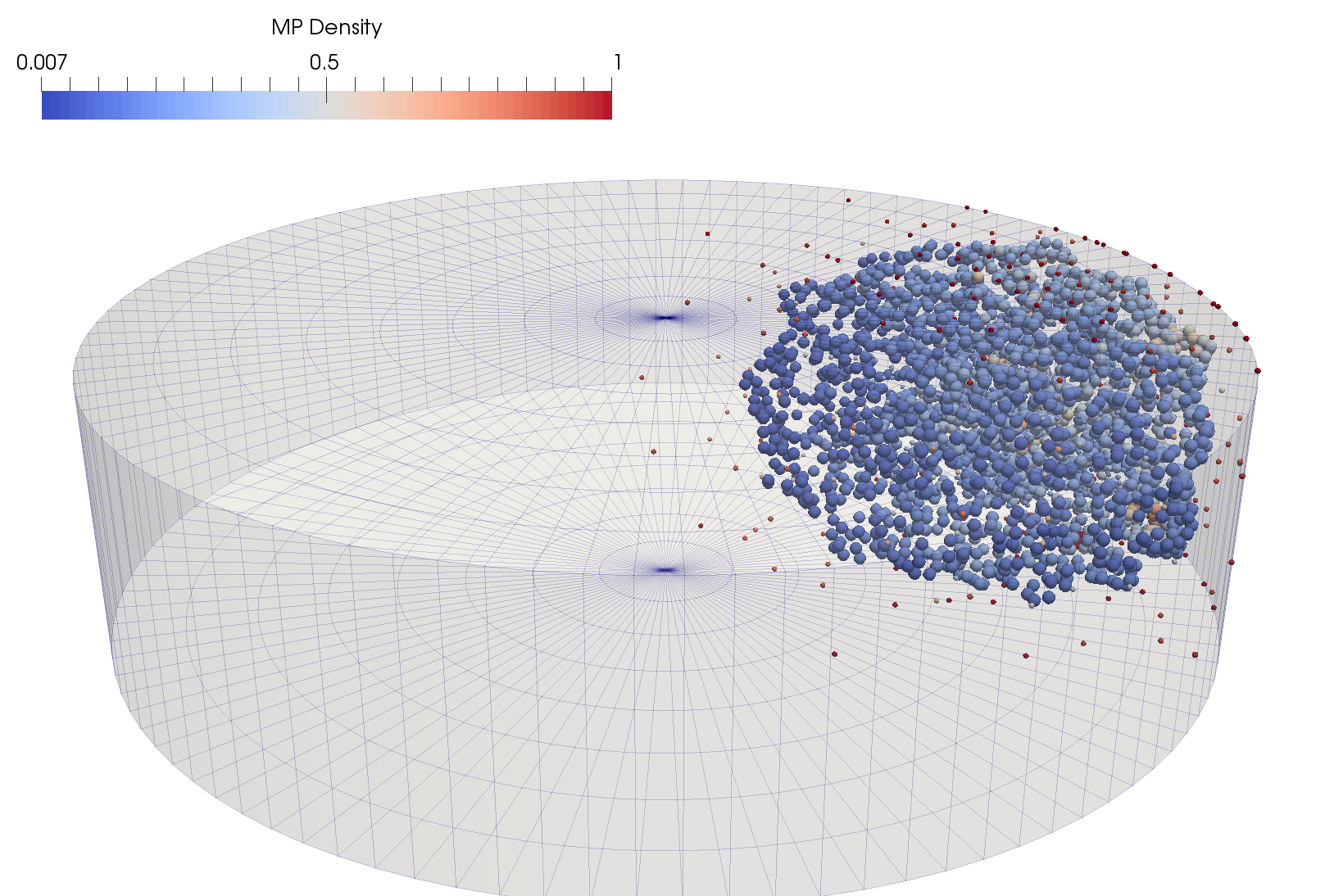}
    \label{fig:b2}}
	\subfigure[$t=0.2$]
    {\includegraphics[width=0.3\textwidth]{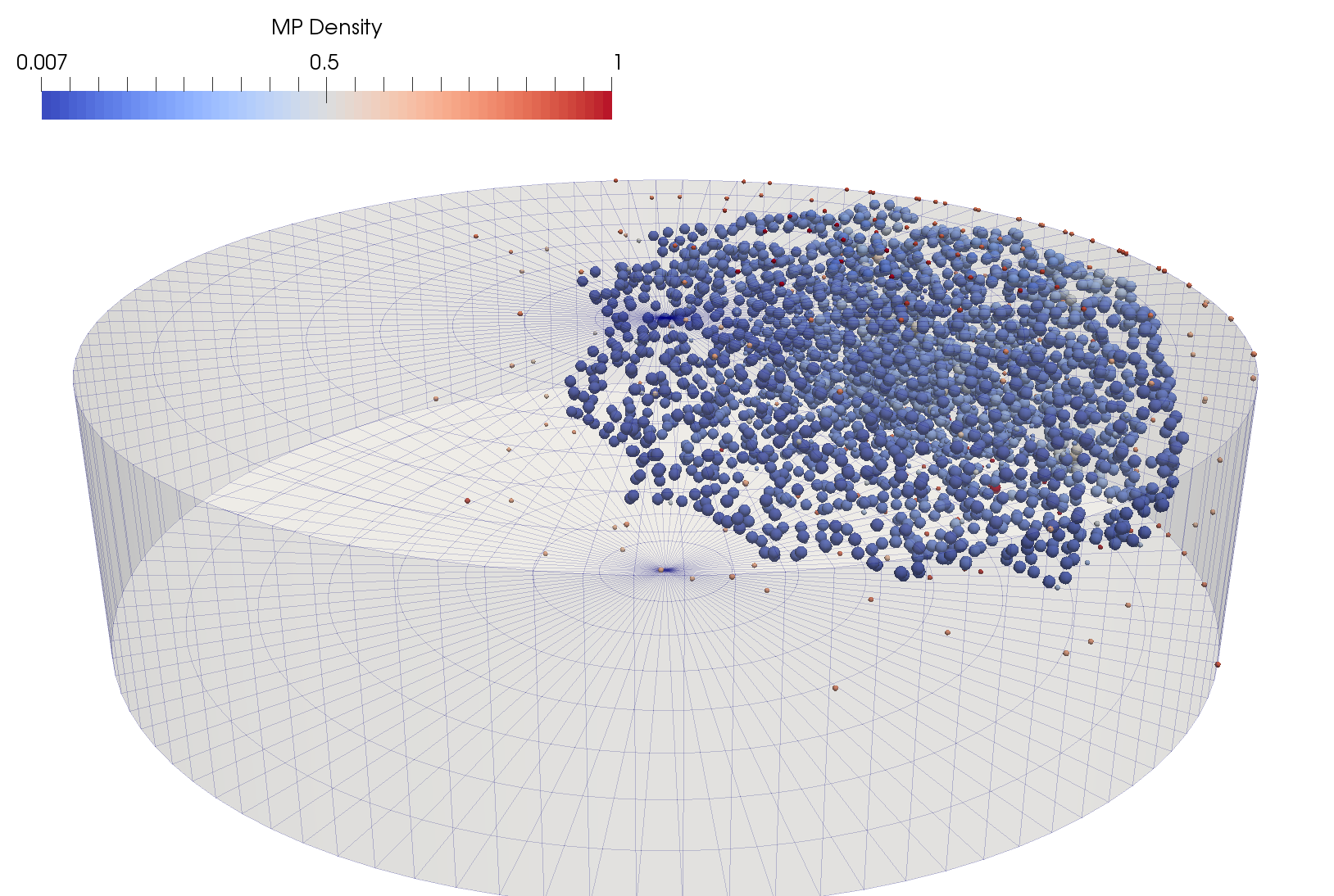}
    \label{fig:b3}}
	\subfigure[$t=0.3$]
    {\includegraphics[width=0.3\textwidth]{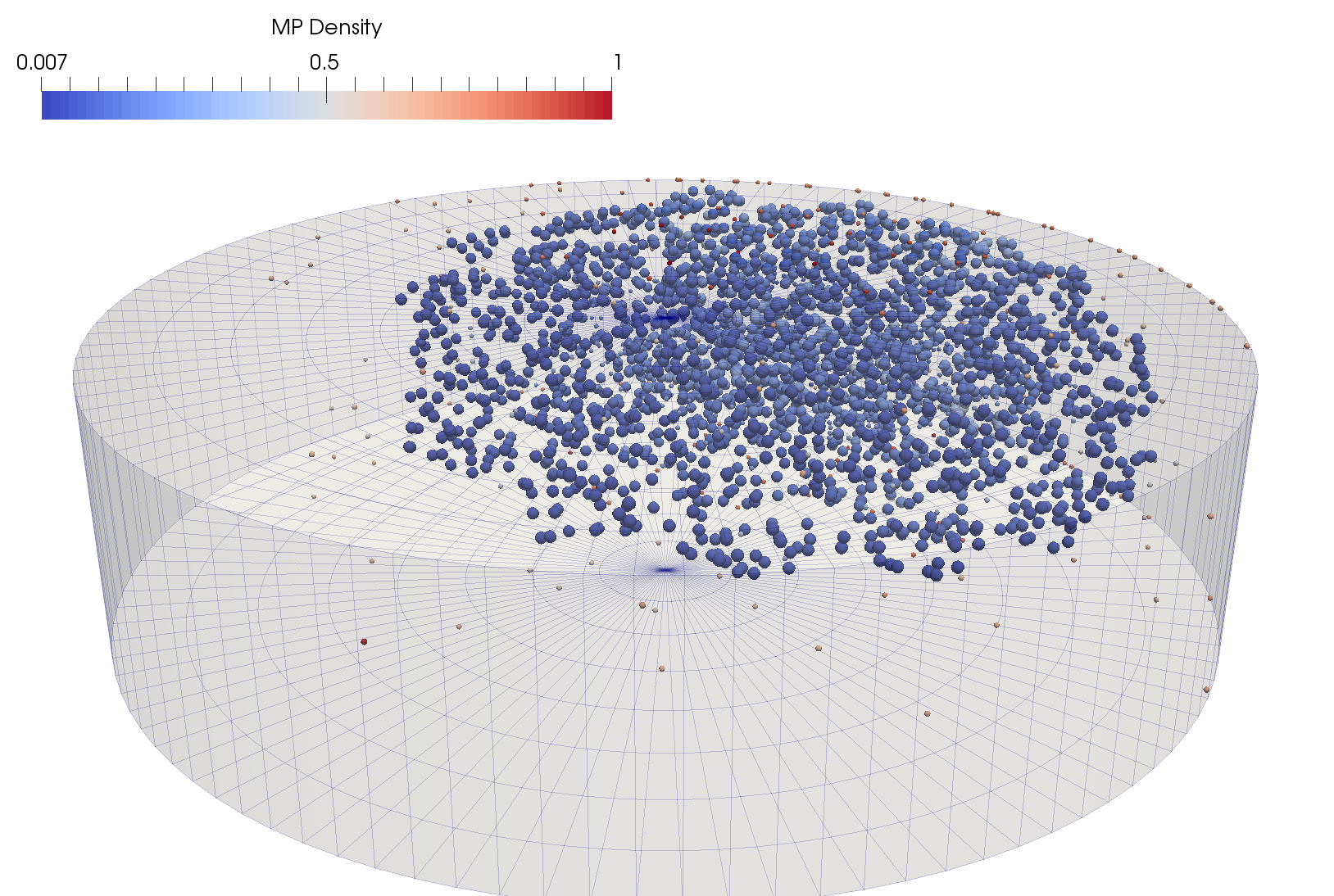}
    \label{fig:b4}}
	\subfigure[$t=1.0$]
    {\includegraphics[width=0.3\textwidth]{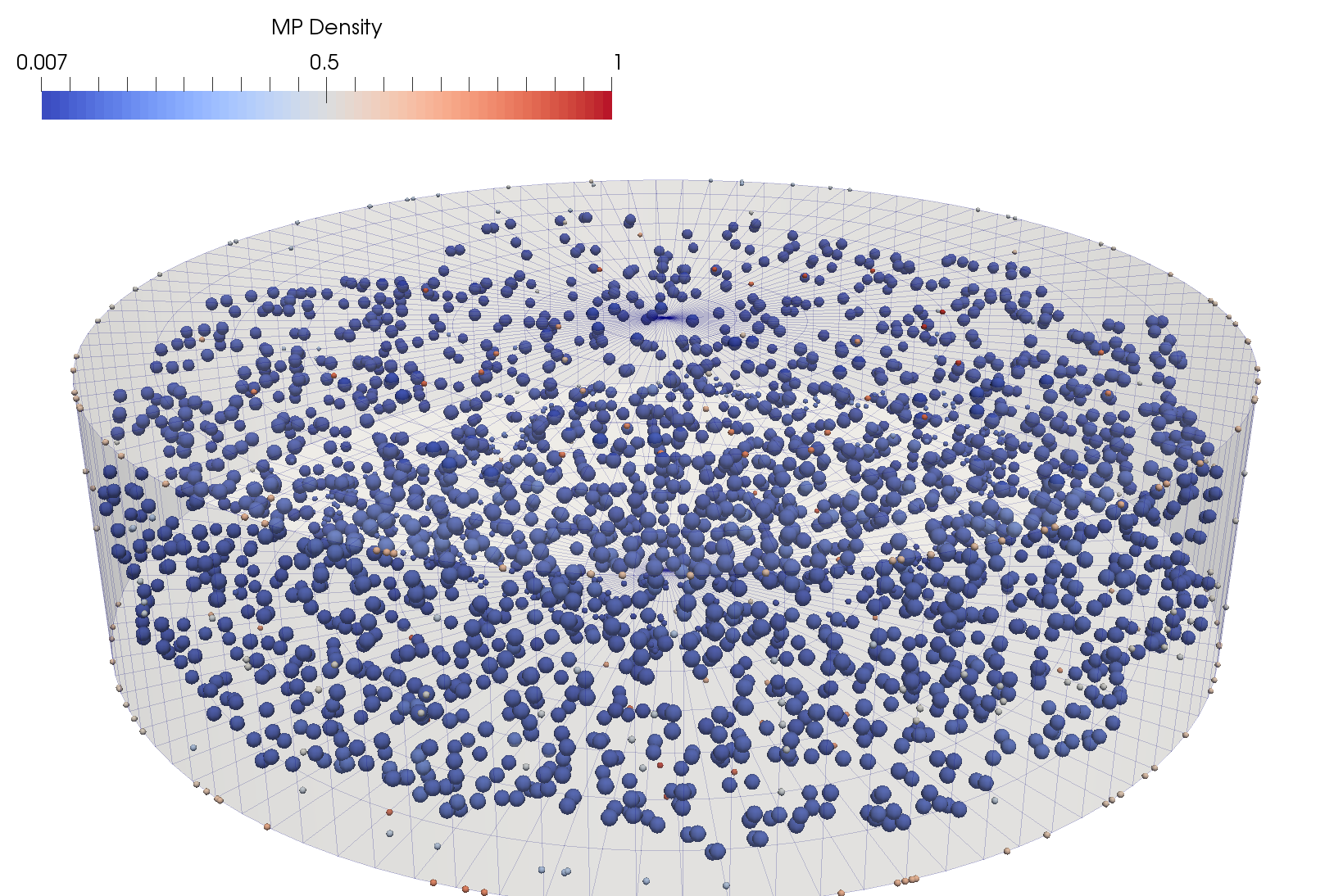}
    \label{fig:b5}}
	\caption{\footnotesize
    Advection-diffusion in a circular cylinder. The initial density is uniform over a spherical region spanning the length of the cylinder. The advection velocity field corresponds to a rigid-body motion at constant angular velocity. Snapshots of the material-point set at three successive times.}
	\label{fig:BucketAdvection}
\end{center}
\end{figure}

The long-term behavior of the method is illustrated by means of the problem of advection-diffusion in a circular cylinder, Fig.~\ref{fig:BucketAdvection}. The radius of the cylinder is $0.5$ and its length is $0.25$. The initial mass density is uniform within a spherical region spanning the length of the cylinder. The density is discretized into 723 nodes and 2,215 material-points. The advection velocity field corresponds to a rigid-body motion at a constant angular velocity of 4. We consider three diffusivities of values $\kappa=0.005$, $\kappa=0.0075$, and $\kappa=0.01$.

Fig.~\ref{fig:BucketAdvection} shows successive snapshots of the material-point distribution for the case of diffusivity $\kappa=0.0075$. As may be observed from the figure, the material and nodal point set spreads from its initial location and tends towards a configuration of uniform density. The long-term stages of the solution are thus characterized by the {\sl fine mixing} of the material points. Again, the robust manner in which the method effects the requisite mixing of material points is noteworthy.

\begin{figure}[htp!]
\begin{center}
	\includegraphics[width=0.45\textwidth]
    {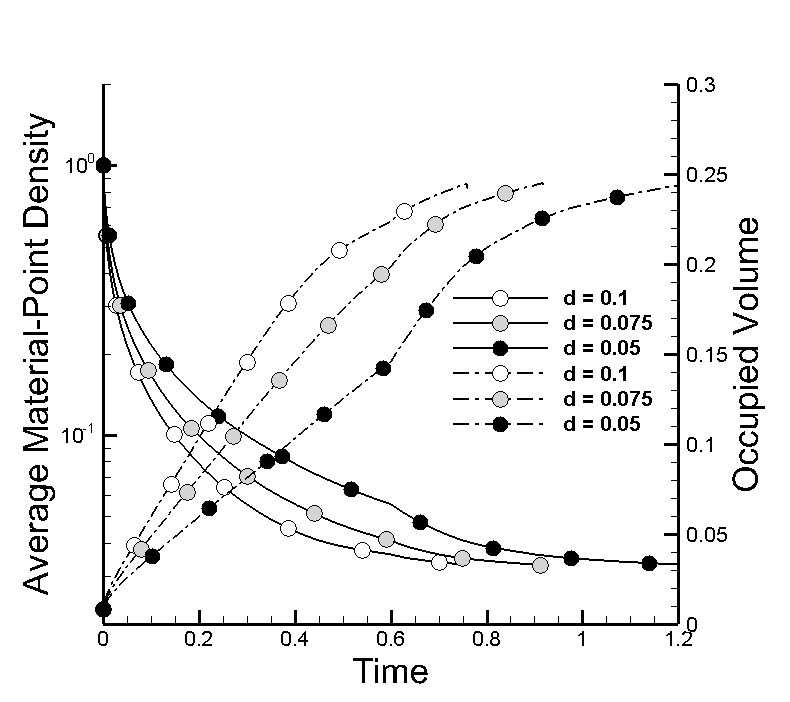}
	\caption{\footnotesize
    Advection-diffusion in a circular cylinder. The initial density is uniform over a spherical region spanning the length of the cylinder. The advection velocity field corresponds to a rigid-body motion at constant angular velocity. Time evolution of the average material-point density and the volume occupied by the material-point set.}
	\label{fig:bucketDensityVolume}
\end{center}
\end{figure}

Fig.~\ref{fig:bucketDensityVolume} shows the time history of the average material-point density and of the volume occupied by the material-point set for the pure advection and the advection-diffusion problems. As expected, the the average density steadily decreases and, correspondingly, the occupied volume increases, with time. However, two distinct regimes are clearly discernible in the evolution of the system. The first stage corresponds to the free spreading of the initial distribution of mass, followed by a second stage corresponding to the mixing of the material points. Interestingly, mixing is observed to accelerate the diffusion process during the mixing phase of the solution.

\section{Concluding remarks}
\label{sec:Conclusions}

We have developed an Optimal Transportation Meshfree (OTM) particle method for advection-diffusion problems regarded as optimal transport of measures. Thus, in sharp contrast to traditional methods of approximation, which regard advection-diffusion as an evolution of a density function taking place is a suitable linear space, here mass distribution is viewed as a {\sl measure} evolving in a non-linear space, or manifold. The method follows after the OTM formulation of \cite{LiHabbalOrtiz2010} for fluid and solid flows and represents a scalar version of that method. As in the vector OTM method, we resort to the incremental variational principle of Jordan, Kinderlehrer and Otto, \cite{JordanKinderlehrerOtto1997, JordanKinderlehrerOtto1998, JordanKinderlehrerOtto1999}, for purposes of time discretization of the diffusive step. This principle characterizes the evolution of the density as a competition between the Wasserstein distance between two consecutive densities and entropy. By virtue of the structure of the Euler-Lagrange equations, sufficiently weakened so as to be linear and undifferentiated in the density, we may approximate the density as a collection of Dirac masses. The interpolation of the incremental transport map, which needs to be conforming, is effected through mesh-free max-ent interpolation \cite{ArroyoOrtiz2006}. Remarkably, the resulting update is {\sl geometrically exact} with respect to advection and volume, which overcomes the chronic problems of spurious diffusion and oscillations that afflict linear-space methods. We present three-dimensional examples of application that illustrate the scope and robustness of the method.

{\black The addition of distributed sources and sinks is straightforward and it simply requires allowing the mass of each particle to change with time. Non-zero Neumann boundary conditions can also be implemented simply by injecting particles or allowing particles to escape through the boundary. Finally, Dirichlet boundary conditions can be enforce by placing fixed particles consistent with a prescribed density on the boundary.}

The accuracy and robustness of the OTM method have been assessed by means of selected advection-diffusion test problems. The numerical tests exhibit clearly the geometrically exact nature of the method with respect to advection and volume. Thus, advection is handled directly by means of a pushforward of the mass density by the advection flow. In addition, the volume of the material points is updated by means of the Jacobian of the incremental transport map. Both operations are geometrically exact, in contrast to the conventional Eulerian treatment of advection and volume through upwinding and the solution of a Poisson problem. The numerical solutions generally exhibit a remarkable robustness and accuracy, even for coarse approximations, enabling, for instance, the simulation of mixing of phases and the propagation of density fronts.

We have noted that the method, as presented here, is subject to spatial and temporal resolution requirements, in particular as regards the need to resolve concave regions of the boundary. However, it should be carefully noted that said resolution requirements are directly tied to the particular implementation adopted here and are not inherent to the optimal transport treatment of advection-diffusion problems. Thus, the convex-programming nature of max-ent interpolation, and its relation to convex domains, may be relaxed by considering signed shape functions \cite{Bompadre2012b}. In addition, the explicit character of the implementation presented here can be generalized by considering incremental minimum principles, other than JKO (\ref{eq:TD:F}), in the spirit of variational integrators \cite{Kane1999, Marsden01, Lew2003, LewMarsdenOrtizWest2004}. This extension enables, in particular, the formulation of implicit and asynchronous time-integration schemes.

Finally, we note the connection between transport problems and gradient flows Wasserstein spaces (e.~g., \cite{LisiniMatthesSavare2012, CavalliNaldi2010}). This connection suggests the applicability of the present approach to broad classes of problems in mechanics and physics, including flow in porous media, Oswald ripening, the Cahn-Hilliard equation, and others. These generalizations and extensions define worthwhile directions for further study but are beyond the scope of this paper.

\ack

LF and MO gratefully acknowledge the support of the U.S.~National Science Foundation through the Partnership for International Research and Education (PIRE) on Science at the Triple Point Between Mathematics, Mechanics and Materials Science, Award Number 0967140.

\begin{appendix}

\section{Elements of optimal transportation theory}

For completeness, in this appendix we collect basic elements of optimal transportation theory that underpin the Optimal Transportation Meshfree (OTM) method presented in the foregoing. A complete account of the theory may be found, e.~g., in \cite{Villani2003, Daneri2010}.

Let $\Omega \subseteq \mathbb{R}^d$ an open convex set and let $\mathscr{P}(\Omega)$, the space of Borel probability measures on $\Omega$. If $\mu \in \mathscr{P}(\Omega)$ and $\varphi: \Omega \rightarrow \Omega$ is a Borel map, $\varphi_{\#}\mu \in \mathscr{P}(\Omega)$ denotes the \emph{push-forward} of $\mu$ through $\varphi$, defined by
\begin{equation}
    \varphi_{\#}\mu(B)
    :=
    \mu(\varphi^{-1}(B)) \quad \text{for every } B \text{ Borel }
    \in \mathscr{P}(\Omega).
\end{equation}
Let $\pi^i, \,\, i = 1 \dots n$ be the canonical projection operator from a product space $\Omega \times \dots \times \Omega$ onto $\Omega$,
\begin{equation}
    \pi^i(z_1,\dots,z_n) := z_i .
\end{equation}
Given $\mu_1 \in \mathscr{P}(\Omega)$ and $\mu_2 \in \mathscr{P}(\Omega)$ the class $\Gamma(\mu_1,\mu_2)$ of \emph{transport plans} or \emph{couplings} between $\mu_1$ and $\mu_2$ is defined as
\begin{equation}
    \Gamma(\mu_1,\mu_2) := \{ \bm{\gamma} \in \mathscr{P}(\Omega \times \Omega) \}: \pi^1_{\#}\bm{\gamma} = \mu_1, \pi^2_{\#}\bm{\gamma} = \mu_2\}.
\end{equation}
We denote by $\mathscr{P}_2(\Omega)$ the space of Borel probability measures with finite second moment: $\mu \in \mathscr{P}(\Omega)$ belongs to $\mathscr{P}_2(\Omega)$ iff:
\begin{equation}
    \int_{\Omega}| x - x_{0}|^2 {\rm d}\mu(x)
    <
    +\infty, \quad \text{for some (and thus every) point }
    x_{0} \in \Omega.
\end{equation}
For every pair of measures $\mu,\nu \in \mathscr{P}_2(\Omega)$ we consider the transport cost
\begin{equation}\label{wdist}
    d_W^2(\mu,\nu)
    :=
    \min \left \{ \int_{\Omega\times \Omega} |x-y|^2 {\rm d} \gamma(x,y) :
    \gamma \in \Gamma(\mu,\nu)\right \}.
\end{equation}
By the direct method of calculus of variations, the minimum problem (\ref{wdist}) admits at least a solution. The minimum value $d_W(\mu,\nu)$ defines a distance between the measures $\mu,\nu \in \mathscr{P}_2(\Omega)$ and the metric space $(\mathscr{P}_2(\Omega),d_W)$ is referred to as the ($L^2$-) \emph{Wasserstein space} on $\Omega$. We denote by $\mathscr{P}_2^a(\Omega)$ the subset of $\mathscr{P}_2(\Omega)$ formed by the absolutely continuous measures with respect to the Lebesgue measure,
\begin{equation}
    \mathscr{P}_2^a(\Omega)
    :=
    \{ \mu \in \mathscr{P}_2(\Omega) : \mu  \ll \mathcal{L}^d \}.
\end{equation}
The following result establishes the existence and uniqueness of optimal transport plans induced by maps, or \emph{optimal transport maps}, got the case in which the initial measure $\mu$ belongs to $\mathscr{P}_2^a(\Omega)$.

\begin{theorem}[Existence and uniqueness of optimal transport maps \cite{Brenier1991, KnottSmith1984} ] For any $\mu \in \mathscr{P}_2^a(\Omega)$ and $\nu \in \mathscr{P}_2(\Omega)$, the Kantorovich optimal transport problem (\ref{wdist}) has a unique solution $\gamma$, which is concentrated on the graph of a transport map, namely, the unique minimizer of Monge's optimal transport problem
\begin{equation}
    \min \left \{ \int_{\Omega} | x - \varphi(x)|^2 {\rm d}\mu(x) \,: \, \varphi_{\#}\mu = \nu  \right \}.
\end{equation}
The map $\varphi$ is cyclically monotone and there exists a convex function $\phi: \Omega \rightarrow \mathbb{R}$ such that $\varphi(x) = \nabla \phi(x)$ for $\mu-a.e. \, x \in \Omega$.
\end{theorem}

We also recall that the manifold $\mu \in \mathscr{P}_2(\Omega)$ can be endowed with a canonical differential structure. In particular, the tangent space to $\mathscr{P}_2(\Omega)$ at the point $\mu$ is defined as
\begin{equation}
    \textnormal{Tan}_{\mu}\mathscr{P}_2(\Omega)
    :=
    \overline{ \{ \xi: \xi = \nabla \eta, \, \eta \in
    C_c^{\infty}(\Omega)\}}^{L^2(\Omega, \mu)} ,
\end{equation}
where the elements of $\textnormal{Tan}_{\mu}\mathscr{P}_2(\Omega)$ may be regarded as velocity fields.  Indeed, for any absolutely continuous curve $t \mapsto \mu_t \in \mathscr{P}_2(\Omega)$ with $\mu_0 = \mu$, there exists a unique $\xi \in \textnormal{Tan}_{\mu}\mathscr{P}_2(\Omega)$ such that
\begin{equation}
    \left. \frac{\partial }{\partial t}\mu_t \right |_{t =0}
    +
    \nabla \cdot (\mu \xi) = 0 \quad \text{in} \,\,\mathcal{D}'(\Omega).
\end{equation}
In particular, if $\mu_t \in \mathscr{P}_2^a(\Omega)$, and $\rho_t$ is its density with $\rho_0 = \rho$, then
\begin{equation}
    \left. \frac{d}{dt} \right |_{t=0} \int_{\Omega} \rho_t(x)\eta(x) dx
    =
    \int_{\Omega} \xi(x) \cdot \nabla \eta(x) \rho(x)dx ,
\end{equation}
for all test functions $\eta \in \mathcal{D}(\Omega)$. Thus, changes in densities in $\mathscr{P}_2^a(\Omega)$ are to be regarded as mass fluxes induced by velocity fields.

We also recall the following differentiation formula. Let $b > a$, $\rho_a$, $\rho_b \in \mathscr{P}_2^a(\Omega)$ and let $\rho_t \in \mathscr{P}_2^a(\Omega)$, $t \in [a,b]$, be such $\rho(\cdot,a) = \rho_a$, $\rho(\cdot,b) = \rho_b$ and
\begin{equation}
    \frac{\partial \rho_t}{\partial t} + \nabla \cdot (\rho_t \xi) = 0 ,
\end{equation}
for some $C^1$, globally bounded, velocity field $\xi$. Then \cite{Villani2003},
\begin{equation}\label{diffeWdist}
    \frac{{\rm d}}{{\rm d}t} d_W^2(\rho_a,\rho_t) \restriction_{t=b}
    =
    2 \int_{\Omega} \langle x - \varphi_{\rho_b},
    \xi(x)\rangle {\rm d}\rho_b(x)
\end{equation}
where $\varphi_{\rho_b}$ is the unique optimal transport map between $\rho_a$ and $\rho_b$.


We conclude by showing that the first variation in $\mathscr{P}_2(\Omega)$ of the entropy
\begin{equation}
    S(\rho_{k+1}) = \int_{\Omega} \rho_{k+1}\log \rho_{k+1} \, dx
\end{equation}
gives (\ref{eq:TD:DFmin}) and (\ref{eq:TD:RVD3}). Thus, let $\xi \in  C^{\infty}_c({\Omega})$ be a smooth vector field with compact support. Define the variation $y = \psi(x,t) = \psi_{t}(x)$ as
\begin{equation}\label{variations}
    \partial_{t} y = \xi(y), \,\, y(0) = x .
\end{equation}
Let $\rho_{t}$ be the pushforward measure of $\rho_{k+1}$ by $\psi_t$, i.~e.,
\begin{equation}\label{push_density}
    \int_{\Omega} \rho_t(x) \eta(x) \, dx
    =
    \int_{\Omega} \rho_{k+1} \eta(\psi_t(x)) \,dx ,
\end{equation}
for all $\eta \in C^{\infty}_{0}(\Omega)$.  Since $\psi_t$ is regular, (\ref{push_density}) is equivalent to the relation between densities
\begin{equation}\label{density_relation}
    \big( {\rm det}(\nabla\psi_t) \, \rho_t \big) \circ \psi_t = \rho_{k+1}.
\end{equation}
In particular,
\begin{equation}
    \int_{\Omega} \rho_t \, dx = \int_{\Omega} \rho_{k+1} \, dx ,
\end{equation}
and $\rho_t$ is in $\mathscr{P}_2(\Omega)$. In addition, we have
\begin{equation}
    \int_{\Omega} \rho_t(x) \log(\rho_t(x)) \, dx 
    {=}
    \int_{\Omega}
    \rho_{k+1}(x)\log(\rho_t(\psi_t(x))) \, dx 
    {=}
    \int_{\Omega} \rho_{k+1} \log \left( \frac{\rho_{k+1}(x)}{\text{det}\nabla \psi_{t}(x)} \right) \, dx ,
\end{equation}
and
\begin{equation}
    \frac{1}{t} \left (  S (\rho_t) - S (\rho_{k+1}) \right )
    =
    -
    \frac{1}{t} \int_{\Omega} \rho_{k+1}(x)\log(\text{det}\nabla \psi_t(x)) \, dx.
\end{equation}
Write
\begin{equation}
    J_{\psi_t} := \text{det}\nabla \psi_t .
\end{equation}
Then, we have
\begin{equation}
\begin{split}
    \frac{{\rm d}}{{\rm d}t}\left [ J_{\psi_t} \right ]|_{t = 0}
    & =
    \left. \left [ \frac{\partial \, J_{\psi_t}}{\partial \, \nabla \psi_t } \,
    {\cdot}
    \, \frac{\partial \, \nabla \psi_t }{\partial t}  \right ] \right |_{t=0}
    =
    \left. \left [ J_{\psi_t} (\nabla \psi_t)^{-T} \, {\cdot} \, \nabla \left(\frac{\partial \psi_t}{\partial t} \right)
    \right ] \right |_{t = 0}
    \\ & =
    \left. \left [ J_{\psi_t} (\nabla \psi_t)^{-T} \, {\cdot} \, \nabla \xi \right ] \right |_{t = 0}
    =
    \text{I} \, {\cdot} \, \nabla \xi = \nabla \cdot \xi ,
\end{split}
\end{equation}
where we have used that $\psi(x,0) = \text{id}$. Taking the limit $t \rightarrow 0$, we find
\begin{equation}
    \frac{{\rm d}}{{\rm d}t} S (\rho_t)|_{t = 0}
    =
    -\int_{\Omega}\rho_{k+1}(x)\nabla \cdot \xi(x) \, dx,
\end{equation}
as advertised.

\end{appendix}


\end{document}